\documentclass{IEEEtran}
\usepackage{cite}
\usepackage{amsfonts}
\usepackage{textcomp}

\usepackage{graphicx, amsmath, theorem, amssymb}
\usepackage{tensor}

\usepackage{algorithm}
\usepackage{algorithmic}

\usepackage{color}

\usepackage{accents}

\newtheorem{theorem}{Theorem}

\newtheorem{lm}{Lemma}

\newtheorem{df}{Definition}

\newtheorem{AS}{Assumption}

\newcommand{\proof}{{\bf Proof:~}}

\newcommand*{\dt}[1]{%
  \accentset{\mbox{\large\bfseries .}}{#1}}

\newcommand{\Li}{\mathcal{L}}
\newcommand{\I}{\mathcal{I}}
\newcommand{\M}{\mathcal{M}}
\newcommand{\E}{\mathcal{E}}
\newcommand{\R}{\mathbb{R}}

\newtheorem{rem}{Remark}

\newsavebox{\mybox}

\def\BibTeX{{\rm B\kern-.05em{\sc i\kern-.025em b}\kern-.08em
    T\kern-.1667em\lower.7ex\hbox{E}\kern-.125emX}}

\begin{document}

\newcommand{\SUIO}{$\mathcal{SUIO}$}
\newcommand{\Obs}{$\mathcal{O}$}
\newcommand{\OBS}{\mathcal{O}}

\newcommand{\Sdi}{\mathcal{UIE}_{d\infty}}
\newcommand{\STk}{\mathcal{UIE}_{m_wk}}
\newcommand{\Sdimk}{\mathcal{UIE}_{mk,d\infty}}
\newcommand{\G}{\mathcal{G}}
\newcommand{\F}{\mathcal{F}}

\newcommand{\xdi}{x_{d\infty}}
\newcommand{\wdi}{w_{d\infty}}
\newcommand{\dentro}{d}
\newcommand{\uideg}{\textnormal{deg}_w}

\newcommand{\Hf}{\mathcal{H}}
\newcommand{\RM}{\mathcal{RM}}
\newcommand{\VS}{$\mathcal{V}$}
\newcommand{\IMU}{$\mathcal{I}$}
\newcommand{\tobs}{\widetilde{\OBS}}
\newcommand{\NO}{s}
\newcommand{\ND}{r}
\newcommand{\tih}{\widetilde{h}}
\newcommand{\C}{\mathcal{C}}
\newcommand{\Mn}{\mathbb{R}^n}

\title{General analytical condition to nonlinear identifiability and its application in viral dynamics}

\author{Agostino Martinelli
\thanks{A. Martinelli is with INRIA Rhone Alpes,
Montbonnot, France e-mail: {\tt agostino.martinelli@inria.fr}} }

\maketitle

\begin{abstract}
Identifiability describes the possibility of determining the values of the unknown parameters that characterize a dynamic system from the knowledge of its inputs and outputs.
This paper finds the general analytical condition that fully characterizes this property.
The condition can be applied to any system, regardless of its complexity and type of nonlinearity. In the presence of time varying parameters, it is only required that their time dependence be analytical.
In addition, its implementation requires no inventiveness from the user as it simply needs to follow the steps of a systematic procedure that 
only requires to perform the calculation of derivatives and matrix ranks.
Time varying parameters are treated as unknown inputs and their identifiability is based on the very recent analytical solution of the unknown input observability problem \cite{SIAMbook,IF22}.
Finally, when a parameter is unidentifiable,
the paper also provides an analytical method to determine infinitely many values for this parameter that are indistinguishable from its true value.
The condition is used to study the identifiability of two nonlinear models in the field of viral dynamics (HIV and Covid-19). 
In particular, regarding the former, a very popular HIV ODE model is investigated,  and the condition allows us to automatically find a  
new fundamental result that highlights a serious error of the current state of the art.

\end{abstract}


\section{Introduction}\label{SectionIntroduction}

Ordinary Differential Equations (ODEs) are used to model a plethora of phenomena in many scientific domains, ranging from the natural and applied sciences up to the social sciences.

An ODE model is characterized by a state that consists of several scalar quantities. Its time evolution is precisely described by a set of ordinary differential equations (one per each state component).
In addition to the state components, an ODE model can also include a set of parameters which may be known or unknown.
Finally, an ODE model is characterized by a set of outputs and, very often, by a set of inputs. The inputs are functions of time that act on the state dynamics and that can be assigned (with some restrictions that depend on the specific case).
The outputs are available quantities and can be expressed in terms of the state and, in some cases, also in terms of the above inputs and the above model parameters.

{\it State observability} and {\it parameter identifiability} are two structural properties of an ODE model. The former characterizes the possibility of inferring the values that the state components take, starting from the knowledge of the system inputs and outputs
(e.g., \cite{Kalman61,Kalman63,Her77,Casti82,Isi95}). The latter characterizes the possibility of inferring the values of the model parameters, again, from the knowledge of the system inputs and outputs
(e.g., \cite{Grew76,Tuna87,Wal96}).

An ODE model may also include the presence of {\it unknown} inputs, i.e., inputs that act on the system dynamics precisely as the aforementioned inputs. However, they differ from them because they are unknown. From a mathematical point of view, an unknown input plays exactly the same role of an unknown time varying parameter. In the literature, the possibility of inferring their values has been called in different manners:
{\it input observability}, {\it input reconstruction}, and sometimes {\it system invertibility}
(e.g., \cite{Basile69,Guido71,Hou98,Krya93,Basile73,Hirs79,Sing81,Fli86}).

Throughout this paper, we adopt the following terminology. {\it Observability} refers to the components of the state. {\it Unknown Input Observability} (abbreviated UIO) 
still refers only to the state components, but when the ODE model also includes the presence of unknown inputs (or time varying parameters).
{\it Identifiability} refers to all the unknown parameters. 
When a parameter is time varying we also use the term {\it unknown input reconstruction} to mean its identifiability.


ODE identifiability analysis is present in a large variety of scientific domains and several interesting methods have been proposed. An exhaustive review of the methods proposed up to 2011 can be found in \cite{Miao11}. 
Some of these methods are very interesting and, in many cases, they can successfully be used to detect all the identifiability properties of a given ODE model. However,
they have the following fundamental limitations:

\begin{itemize}

\item They are not general. They are based on several properties that sometimes cannot be exploited in the presence of a given type of system nonlinearity (i.e., a given function $f$ and/or $h$ in Equation (\ref{EquationSystemDefinitionUIOGeneral})). In particular, many of them only hold for a polynomial nonlinearity.

\item They cannot be executed automatically. In most of cases, they must be adapted to the specific case under investigation and this adjustment requires inventiveness by the user. In particular, this process cannot be carried out by simply following the steps of a systematic procedure (e.g., by running a code without human intervention).

\end{itemize}

The second limitation becomes particularly relevant in the case of time varying parameters where, in some cases, an erroneous application of these methods provided wrong results
(in our study about the identifiability of the ODE HIV model given in Section \ref{SectionHIV}, we obtain results that contradict the ones available in the state of the art, which are obtained by using one of these methods).

In the last decade, 
thanks to the progress in UIO, several new interesting approaches have been proposed.
Specifically,
the nonlinear UIO problem was
approached by introducing an extended state that includes the original state together with the unknown inputs and their time derivatives up to a given order \cite{Belo10,MED15,Villa19b,Maes19,Villa16b}. 
In \cite{Villa19b,Maes19,Villa16b}, this was used to study the identifiability of the unknown time varying parameters. 
In particular, in \cite{Villa19b,Maes19,Villa16b}, several
automatic iterative algorithms were introduced. These algorithms work automatically and are able to check the identifiability of the time varying parameters. 
On the other hand, they suffer from the following  limitations:

\begin{itemize}

\item They do not converge, automatically. 
In particular, at each iterative step, the state is extended by including new time derivatives of the unknown inputs. Consequently, if the extended state is observable at a given step, convergence is achieved. However, if this were not the case, we can never exclude that, at a later step, the extended state becomes observable. Therefore, in the presence of unobservability, all these algorithms remain inconclusive.

\item Due to the previous state augmentation, the computational burden can easily become prohibitive after a few steps.

\end{itemize}

Very recently, we introduced the analytical solution of the nonlinear UIO \cite{IF22}, starting from the results/derivations presented in \cite{SIAMbook}\footnote{This analytical solution is also based on the results obtained in \cite{TAC19} and \cite{SARAFRAZI}, which deal with special systems, namely characterized by a single unknown input, and dynamics linear with respect to this unknown input.}.
In particular, the systematic procedure introduced in \cite{IF22} does not encounter the aforementioned limitations\footnote{Note that, the solution introduced in \cite{IF22} could need the inclusion of some of the UIs in the state. However, this inclusion is targeted and terminates in a finite number of steps, as soon as the extended system achieves its {\it highest unknown input degree of reconstructability} (see Section 5 in \cite{IF22}).} and provides a full answer to the problem of the state observability for any nonlinear ODE model, in the presence of unknown inputs (or time varying parameters). On the other hand, this systematic procedure deals with the observability of the state and not with the identifiability of the time varying parameters.
In \cite{IF22}, we exploited this analytical solution to provide a preliminary result also about the identifiability of the time varying parameters (unknown input reconstruction).
One of the goals of this paper is the extension of this preliminary and partial result, and the introduction of the general analytical condition to study the identifiability of the time varying parameters. Then, we adopt this condition to study the identifiability of a very popular HIV model (e.g., \cite{Miao11,Villa19b,Per99,Per02,Chen08,Moog10,Xia03}) and a Covid-19 model \cite{VillaCovid,SEIARPribylova}.


The contributions of this paper are the following three:

\begin{enumerate}

\item Introduction of the analytical condition that fully characterizes the local identifiability of all the unknown parameters of any ODE model (Section \ref{SectionCondition}).

\item Determination of the continuous transformations (one parameter Lie groups) that allow us to determine infinitely many values, for any unidentifiable parameter, that are indistinguishable from its true value.

\item 
Determination of the identifiability properties of a very popular HIV model and a Covid-19 model, with the detection of a serious error in the state of the art on the identifiability of the HIV model.

\end{enumerate}

Regarding the first contribution, the time varying parameters are divided in two distinct groups. The second group can be empty. The parameters that belong to it are always not identifiable, even locally (first part of Theorem \ref{TheoremIdentifiabilityCan}). The condition  that fully characterizes the local identifiability of the parameters of the first group is given in the first part of Theorem \ref{TheoremIdentifiabilityGeneral} (Equation (\ref{EquationConditionGeneral})).
Finally, the condition that characterizes the local identifiability of the constant parameters is given in Section \ref{SubSectionConditionConstant} (Equation (\ref{EquationConditionConstant})).

Regarding the second contribution, the continuous transformations are given by the system of differential equations in (\ref{EquationConditionDiffEqSystem}) for the time varying parameters of the first group and by (\ref{EquationConditionTransformationLast}) for the time varying parameters of the second group. Note that, the system of differential equations in (\ref{EquationConditionDiffEqSystem}) also provides the continuous transformations for the constant parameters (see Remark \ref{RemarkContinuousConstant} in Section \ref{SubSectionConditionGeneral}).

Regarding the third contribution,
our results improve previous results in the state of the art. In particular, in both cases, the condition in (\ref{EquationConditionGeneral}) provides the unidentifiability of the time varying parameter.
In other words, in both cases, the time varying parameter cannot be uniquely identified (for the HIV model this result is in contrast with the state of the art).
~In addition, for these ODE models it is possible to obtain an analytical solution of the system of differential equations in (\ref{EquationConditionDiffEqSystem}).
This allows us to determine infinitely many values of the unidentifiable parameter that agree with the same outputs of the model. 
Regarding the HIV ODE model, in Section \ref{SubSectionHIVComparisonSOTA},  we directly prove that the solutions of (\ref{EquationConditionDiffEqSystem}) produce the same outputs. This unequivocally shows the unidentifiability of the system, in contrast with the results available in the state of the art.
Finally, we determine the minimal external information (external to the knowledge of the outputs) requested to uniquely determine the system parameters.

The paper is organized as follows.
Section \ref{SectionBasic} provides a basic mathematical characterization of the problem together with the assumptions.
Section \ref{SectionCondition} provides the analytical condition (all the proofs are given in the appendix).
Sections \ref{SectionHIV}, and \ref{SectionCovid} investigate the very popular HIV model discussed in \cite{Per99,Per02,Xia03,Chen08,Moog10,Miao11,Villa19b}, and
the Covid-19 model introduced in \cite{VillaCovid,SEIARPribylova}. 
Section \ref{SectionVISFM} provides an application in robotics to further illustrate the power and the generality of the first two paper contributions.
Specifically, we study the identifiability of a very popular perception model, which is based on visual and inertial sensing. 
~Finally, Section \ref{SectionConclusion} provides our conclusion.


\section{Problem and basic assumptions}\label{SectionBasic}

We start from the very general model:

\begin{equation}\label{EquationSystemDefinitionUIOGeneral}
\left\{\begin{array}{ll}
 \dot{X} &=  f(X(t), ~t, ~U(t), ~Q, ~W(t))\\
  y &= h(X(t), ~t, ~U(t), ~Q, ~W(t)), \\
\end{array}\right.
\end{equation}

where:

\begin{itemize}

\item The functions $f$ and $h$ are defined on $\M_m\times\I\times\M_{m_u}\times\M_q\times\M_{m_w}$, with $\M_m$, $\M_{m_u}$, $\M_q$, and $\M_{m_w}$,  $\C^\infty-$manifolds of dimension $m$, $m_u$, $q$, and $m_w$, respectively, and $\I\subseteq\mathbb{R}$  an open time interval.

\item $X(t)\in\M_m$ is the state. 

\item $y(t)=[y_1(t),~\ldots,~y_p(t)]^T$ is the output vector.

\item $U(t)=\left[U_1(t),~\ldots,~U_{m_u}(t)\right]^T\in\M_{m_u}$ is the known input vector.

\item $Q=[Q_1,~\ldots,~Q_q]^T\in\M_q$ is the set of the unknown constant parameters.

\item $W(t)=\left[W_1(t),~\ldots,~W_{m_w}(t)\right]^T\in\M_{m_w}$ is the set of the unknown time varying parameters.

\end{itemize}
Note that $W(t)$ can be regarded as a system input vector that, precisely as $U(t)$, acts on the system dynamics. However, it differs from $U(t)$ in two fundamental respects: (i) it cannot be assigned, and (ii) it is unknown.

We make the following assumptions:

\begin{AS}
\label{Assumptionfh}
The functions $f$ and $h$ are smooth in all their arguments.
\end{AS}

\begin{AS}
\label{AssumptionW}
The function $W(t)$ is analytic in $t$.
\end{AS}

\begin{AS}[Only needed if $\boldsymbol{f}$ is not affine in $\boldsymbol{U}$]
\label{AssumptionU}
The function $U(t)$ is differentiable.
\end{AS}


The problem that we solve in this paper is the introduction of the analytic condition that fully characterizes the local identifiability of the system parameters, both constant and time varying. 
This is obtained by exploiting recent results on the unknown input observability problem.
To exploit these results, instead of (\ref{EquationSystemDefinitionUIOGeneral}) we adopt the following system characterization:

\begin{equation}\label{EquationSystemDefinitionUIO}
\left\{\begin{array}{ll}
  \dot{x} &=   g^0(x, t)+\sum_{k=1}^{m_u}f^k (x, t) u_k(t) +  \sum_{j=1}^{m_w}g^j (x, t) w_j(t)  \\
  y &= [h_1(x, t),\ldots,h_p(x, t)], \\
\end{array}\right.
\end{equation}

where:

\begin{itemize}

\item $x\in\M_n$ is the state and $\M_n$ is a $\C^\infty-$manifold of dimension $n$. 


\item $u_1(t), \ldots,u_{m_u}(t)$ are the known inputs.

\item $w_1(t), \ldots,w_{m_w}(t)$ are the unknown inputs or the unknown time varying parameters.


\item $f^1,\ldots,f^{m_u},g^0,g^1,\ldots,g^{m_w}$ are smooth vector functions of $x$ and $t$.


\item $h_1,\ldots,h_p$ are smooth scalar functions of $x$ and $t$.

\end{itemize}

Note that, for our purposes, the above characterization can easily account for the presence of constant parameters ($Q$) by including all of them in the state $x$ and by suitably setting $f^1,\ldots,f^{m_u},g^0,g^1,\ldots,g^{m_w}$, namely, by setting to zero all their components that will give 
 $\dt{Q}$ in (\ref{EquationSystemDefinitionUIO}).
Note that many nonlinear systems have the structure in (\ref{EquationSystemDefinitionUIO}). When this is not directly the case, it is possible to easily convert (\ref{EquationSystemDefinitionUIOGeneral}) to (\ref{EquationSystemDefinitionUIO}).
~We set $u(t)=\dt{U}$, and $w(t)=\dt{W}$ (with $u=[u_1,\ldots,u_{m_u}]$ and $w=[w_1,\ldots,w_{m_w}]$) and we include $U(t)$, $W(t)$ and $Q$ in the state (i.e., $x=[X^T, ~U^T, ~W^T, ~Q^T]^T$). 
In the rest of this paper, we directly refer to the characterization given in (\ref{EquationSystemDefinitionUIO}).

\vskip.2cm
Given the system in (\ref{EquationSystemDefinitionUIO}) with Assumptions \ref{Assumptionfh} and \ref{AssumptionW}
(which is equivalent to the system in (\ref{EquationSystemDefinitionUIOGeneral}) when also Assumption \ref{AssumptionU} is honoured), the goal of this paper is to
provide the analytical condition to check the local identifiability of all its unknown parameters.
%
%
%

\section{General analytical condition for nonlinear identifiability}\label{SectionCondition}

To provide the analytical condition that characterizes the local identifiability, we must first execute the systematic procedure introduced in \cite{IF22}, which provides the state observability in the presence of unknown inputs (Algorithm 9 in \cite{IF22}).
For this reason, before introducing this condition, we must remind the reader of some concepts and definitions introduced in \cite{IF22} (Section \ref{SubSectionConditionReminders}) and of some features of Algorithm 9 in \cite{IF22} (Section \ref{SubSectionConditionProcedureUIO}).
Then, Sections \ref{SubSectionConditionWhenOBS} and \ref{SubSectionConditionGeneral} provide the condition that characterizes the local identifiability of the time varying parameters.
Section \ref{SubSectionConditionWhenOBS} deals with the special case when the state is observable (Theorem \ref{TheoremIdentifiabilityObservable}) and Section \ref{SubSectionConditionGeneral} provides the condition in the general case (Theorems \ref{TheoremIdentifiabilityGeneral} and \ref{TheoremIdentifiabilityCan}). 
Note that the results stated by all these theorems hold for a system (denoted by $\E$, in this paper) that may differ from the original system. In particular, during the execution of Algorithm 9 in \cite{IF22},
some of the unknown inputs of the original system could have been included in the state and, consequently, the corresponding unknown inputs of $\E$ become the time derivatives of the original unknown inputs. On the other hand, even when this is the case, it is immediate to obtain the identifiability properties of the original system once we have obtained the identifiability properties of $\E$, as we explain in Section \ref{SubSectionConditionSigma}.
Finally, Section \ref{SubSectionConditionConstant} provides the condition that characterizes the local identifiability of the constant parameters.
All the proofs are given in the appendix.

\subsection{Basic concepts introduced in \cite{IF22}}\label{SubSectionConditionReminders}

We remind the reader of several definitions and concepts introduced in \cite{IF22} (where the reader is addressed for further details).

\subsubsection{Unknown input reconstructability matrix}\label{SubSubSectionConditionUIRecMatrix}
Given a system characterized by (\ref{EquationSystemDefinitionUIO}) and $k$ scalar and smooth functions of the state, $\lambda_1(x),\ldots, \lambda_k(x)$, the unknown input reconstructability matrix of the system from $\lambda_1,\ldots, \lambda_k$ is:

\[
\mathcal{RM}\left( \lambda_1,\ldots, \lambda_k\right)
:=
\left[\begin{array}{cccc}
\Li_{g^1} \lambda_1 & \Li_{g^2} \lambda_1 & \ldots & \Li_{g^{m_w}} \lambda_1 \\
\Li_{g^1} \lambda_2 & \Li_{g^2} \lambda_2 & \ldots & \Li_{g^{m_w}} \lambda_2 \\
\ldots &\ldots &\ldots &\ldots \\
\Li_{g^1} \lambda_k & \Li_{g^2} \lambda_k & \ldots & \Li_{g^{m_w}} \lambda_k \\
\end{array}
\right]
\]
where $\Li_\tau$ is the Lie derivative operator along the vector $\tau$. Note that $\lambda_1,\ldots, \lambda_k$ can also depend on the time $t$.

\subsubsection{Unknown input degree of reconstructability}\label{SubSubSectionConditionUIDeg}
Given a system characterized by (\ref{EquationSystemDefinitionUIO}), and the functions $\lambda_1,\ldots, \lambda_k$, the unknown input degree of reconstructability of the system from $\lambda_1,\ldots, \lambda_k$ is the rank of $\mathcal{RM}\left(\lambda_1,\ldots, \lambda_k\right)$\footnote{As in \cite{IF22}, in this paper we always refer to an open set of $\M_n$ where the above rank takes a constant value.}. When the set of functions $\lambda_1,\ldots, \lambda_k$ consists of all the independent observable functions\footnote{The definition of observable function can be found in \cite{TAC19}. An observable function is a scalar field constant on the indistinguishable sets (see Definition 2 in \cite{TAC19}). From a practical point of view, an observable function is a scalar field such that the values that it takes during a given time interval can be uniquely determined from the system inputs and outputs. For instance, let us suppose that the state that characterizes a vehicle that moves on a plane includes its Cartesian coordinates $x$ and $y$ and let us suppose that, from the inputs and outputs, we cannot determine their value but we can obtain the distance of the vehicle from the origin. The state is unobservable. However, the function $\sqrt{x^2+y^2}$ is an observable function.}$^,$\footnote{Here, with {\it independent} we mean that their differentials are independent covectors. Hence, we can have at most $n$ independent observable functions (when there are $n$ independent observable functions the entire state is evidently observable).} we omit to specify the functions and we refer to the unknown input degree of reconstructability of the system.

\subsubsection{Canonic system with respect to its unknown input and canonical form}\label{SubSubSectionConditionCanonic}
By construction, given a system characterized by (\ref{EquationSystemDefinitionUIO}), its unknown input degree of reconstructability from any set of scalar functions cannot exceed $m_w$.
A system is canonic with respect to its unknown inputs if its unknown input degree of reconstructability is $m_w$. 
In general, at the beginning, the observable functions are not available, with the exception of the system outputs, $h_1,\ldots,h_p$, (which are certainly observable functions).
We say that the system has been set in canonical form with respect to its unknown inputs as soon as the unknown input degree of reconstructability from
all the available observable functions is $m_w$\footnote{When we perform an observability analysis we can somehow get new observable functions, in addition to the outputs. 
~This is precisely the case when we follow the steps of Algorithm 9 in \cite{IF22}.}. 
Note that there are systems that are not canonic with respect to their unknown inputs and, consequently, cannot be set in canonical form.

\subsubsection{Finite unknown input extension}\label{SubSubSectionConditionUIE}
Given a system characterized by (\ref{EquationSystemDefinitionUIO}) we call a finite unknown input extension of it any system that is obtained by including in the state some (or all) its unknown inputs together with the time derivatives up to a given order (the order can be different for each unknown input included). The resulting system, still satisfies (\ref{EquationSystemDefinitionUIO}) with $m_u$ known inputs (which remain the same) and $m_w$ unknown inputs (which will become the time derivatives of a given order of the original ones).
The new dynamics is characterized by new vectors that can be easily obtained from (\ref{EquationSystemDefinitionUIO})  (see Equations (21)-(23) in Section 5 of \cite{IF22}).

\subsubsection{Highest unknown input degree of reconstructability}\label{SubSubSectionConditionHDeg}
Given a system characterized by (\ref{EquationSystemDefinitionUIO}), its highest unknown input degree of reconstructability is the largest unknown input degree of reconstructability of all its finite unknown input extensions. 

The execution of Algorithm 9 in \cite{IF22} can set up a finite unknown input extension of the original system (in many cases it is unnecessary).
When necessary, the state augmentation is minimal, in the sense that only the unknown inputs needed to achieve the highest unknown input degree of reconstructability are included. In addition, the highest unknown input degree of reconstructability is always achieved in a finite number of steps\footnote{This is proved in \cite{IF22} (Section 7.2, Proposition 2, and Lemma 1 in Section 7.1).}. Algorithm 9 in \cite{IF22} automatically selects which unknown inputs must be included.
In this paper, we denote by $\E$ this special unknown input extension.

Note that $\E$ can differ from the original system because some of its unknown inputs have become the time derivatives of a given order of the original unknown inputs. This is the case when the highest unknown input degree of reconstructability of the original system does not coincide with its unknown input degree of reconstructability. In addition, $\E$ may differ from the original system because its unknown inputs could have been re-ordered and, consequently, the vector fields $g^1,\ldots,g^{m_w}$ in (\ref{EquationSystemDefinitionUIO}) exchanged with each other.

\subsection{Outcomes of Algorithm 9 in \cite{IF22}}\label{SubSectionConditionProcedureUIO}
Algorithm 9 in \cite{IF22} is a recursive procedure that provides, in a finite number of steps, the following outcomes:

\begin{enumerate}

\item The aforementioned unknown input extension $\E$.

\item The unknown input degree of reconstructability of $\E$, denoted by $m$, together with a set of observable functions, $\widetilde{h}_1,\ldots,\widetilde{h}_m$,
such that the unknown input reconstructability matrix of $\E$ from $\widetilde{h}_1,\ldots,\widetilde{h}_m$ is full rank.
Note that if $\E$ is canonic with respect to its unknown inputs then $m=m_w$.

\item The observability codistribution of $\E$, denoted by $\OBS$.

\end{enumerate}

\subsection{Identifiability of the time varying parameters when the state is observable}\label{SubSectionConditionWhenOBS}
In \cite{IF22}, we proved
that, when the state that characterizes the system is observable, if the system is canonic with respect to its unknown inputs, then all the unknown inputs can be reconstructed (Theorem 4 in \cite{IF22}).
It also holds the viceversa and, consequently, when the state is observable, we obtain a complete answer to our problem which is given by the following theorem: 


\begin{theorem}\label{TheoremIdentifiabilityObservable}
Let us suppose that the state of $\E$ is observable. $\E$ is canonic with respect to its unknown inputs if and only if all the unknown inputs are locally identifiable.
\end{theorem}

\proof{The proof can be obtained by remarking that the statement is a consequence of the results stated by Theorem \ref{TheoremIdentifiabilityGeneral} and Theorem \ref{TheoremIdentifiabilityCan} (see Remark \ref{Remark2And3Include1}). Due to its simplicity, we also provide a direct proof in Appendix \ref{SubSectionAppendixProofTheoremIdentifiabilityObservable}.
 $\blacktriangleleft$}

Note that, 
~Algorithm 9 in \cite{IF22} provides the highest unknown input degree of reconstructability ($m$). $\E$ is canonic with respect to its unknown input if and only if $m=m_w$.
Therefore, Theorem \ref{TheoremIdentifiabilityObservable}
provides a full answer to the problem of unknown input reconstruction when the state is observable.

%

\subsection{Identifiability of the time varying parameters in the general case}\label{SubSectionConditionGeneral}

Let us denote by $\OBS^\bot$ the orthogonal distribution of $\OBS$. Note that, when the state is observable,
$\OBS^\bot$ only contains the null vector. We remind the reader that $m$ is the unknown input degree of reconstructability of $\E$ (or the highest unknown input degree of reconstructability of the original system).
The following fundamental result fully characterizes the local identifiability of the first $m$ unknown inputs of $\E$: 

\begin{theorem}\label{TheoremIdentifiabilityGeneral}
The unknown input $w_j$ ($j=1,\ldots,m$) is locally identifiable if and only if, for any $\xi\in\OBS^\bot$ and for any $\alpha=0,1,\ldots,m$, we have:

\begin{equation}\label{EquationConditionGeneral}
\sum_{i=1}^m~\nu^i_j \xi^\alpha_i=0,
\end{equation}

where: 

\begin{itemize}

\item $\nu$ is the inverse of the following $(1,1)$-tensor field\footnote{Note that, in \cite{IF22} the tensor $\nu$ was denoted by $^m\nu$ when $m<m_w$. Similarly, the tensor $\mu$ was denoted by $^m\mu$ when $m<m_w$.}:

\begin{equation}\label{EquationConditionMu}
\mu^i_j = \mathcal{L}_{g^i}\widetilde{h}_j, ~~~i,~ j=1,\ldots,m
\end{equation}

\item $\xi^\alpha_i$ depends on $\xi$ and it is, for any $i=1,\ldots,m$:

\begin{equation}\label{EquationConditionCoeffXi}
\xi^\alpha_i=\left\{
\begin{array}{ll}
\left(
\frac{\partial}{\partial x}
\dt{\Li}_{g^0}\widetilde{h}_i
\right)
\cdot\xi,& \alpha=0\\
\left(
\frac{\partial}{\partial x}\Li_{g^\alpha}\widetilde{h}_i
\right)
\cdot\xi,& \alpha=1,\ldots,m\\
\end{array}
\right.
\end{equation}
with $\dt{\Li}_{g^0}=\Li_{g^0}+\frac{\partial}{\partial t}$.

\end{itemize}

In addition, let us denote by $x(t)$ the true state of $\E$ at time $t$ and by $w_1(t),\ldots,w_m(t)$ the first $m$ true unknown inputs of $\E$ at time $t$. For any vector field $\xi\in\OBS^\bot$, let us consider the following system of differential equations in $\tau$ (where $t$ is a fixed parameter):

\begin{equation}\label{EquationConditionDiffEqSystem}
\left\{\begin{array}{ll}
  \frac{dx'}{d\tau} &=   \xi(x'(t,~\tau))\\
  \frac{dw'}{d\tau} &=   ~\chi(x'(t,~\tau),~w'(t,~\tau))\\
  x'(t,~0) &=x(t),~~w'(t,~0)=w(t),\\
\end{array}\right.
\end{equation}
with $\chi=\left[\chi_1,\ldots,\chi_m
\right]^T$ and $\chi_j$, $j=1,\ldots,m$, is:
\begin{equation}\label{EquationConditionChij}
 \chi_j =\chi_j(x,~w) = -\sum_{i=1}^{m}~\nu^i_j \left(\xi^0_i +  \sum_{k=1}^{m}\xi^k_i w_k \right).
\end{equation}
Let us assume that the solution of (\ref{EquationConditionDiffEqSystem}) exists for any $t\in\I$ and for any $\tau\in\mathcal{U}$, with $\mathcal{U}$ an open interval of $\R$ that includes $\tau=0$. Then, the solution of (\ref{EquationConditionDiffEqSystem}) at any $\tau\in\mathcal{U}$, i.e., $x'(t,~\tau)$ and $w'(t,~\tau)$, is indistinguishable from $x(t)$ and $w(t)$.
%
%
\end{theorem}

\proof{The proof is given in Appendix \ref{SubSectionAppendixProofTheoremIdentifiabilityGeneral}.}

When the system is canonic with respect to the unknown inputs, $m=m_w$ and the above result fully characterizes the local identifiability of {\it all} the time varying parameters. On the other hand, when $m<m_w$, Theorem \ref{TheoremIdentifiabilityGeneral}
provides no answer about the identifiability of the last $m_w-m$ time varying parameters. The answer for them is trivial. All of them are not locally identifiable. In particular, we have the following result:

\begin{theorem}\label{TheoremIdentifiabilityCan}
Let us suppose that $m<m_w$ (system not canonic with respect to the unknown inputs). The last $m_w-m$ time varying parameters are not identifiable, even locally.
In addition, for any $\tau\in\R$, $w_j(t)$ is indistinguishable from $w'_j(t,~\tau)$, defined as follows:

\begin{equation}\label{EquationConditionTransformationLast}
w'_j(t,~\tau):=w_j(t)+\tau, ~~~~~j=m+1,\ldots,m_w.
\end{equation}
\end{theorem}

\proof{The proof is given in Appendix \ref{SubSectionAppendixProofTheoremIdentifiabilityCan}.}

\vskip .2cm

\noindent We provide the following remarks.

\begin{rem}\label{RemarkHowToSolve}
The solution of (\ref{EquationConditionDiffEqSystem}) is obtained by first determining $x'(t,~\tau)$ (i.e., by first solving $ \frac{dx'}{d\tau} =   \xi(x'(t,~\tau))$, with initial condition $x'(t,~0) =x(t)$), and then by using it  in $ \frac{dw'}{d\tau} =   ~\chi(x'(t,~\tau),~w'(t,~\tau))$, which becomes a differential equation of only $w'(t,~\tau)$. 
\end{rem}

\begin{rem}\label{RemarkContinuous}
Theorem \ref{TheoremIdentifiabilityGeneral} not only provides a simple condition to check the local identifiability. It also provides a possible manner to quantitatively determine infinitely many values for them, which are indistinguishable from the true ones (the solution of (\ref{EquationConditionDiffEqSystem})). Similarly, Theorem \ref{TheoremIdentifiabilityCan} determines the trivial expression of infinitely many values of the last $m_w-m$ time varying parameters indistinguishable from their true values (Equation (\ref{EquationConditionTransformationLast})).
\end{rem}

\begin{rem}\label{RemarkContinuousConstant}
The solution of (\ref{EquationConditionDiffEqSystem}) determines indistinguishable values for both the first $m$ time varying parameters and the constant parameters. The latter are indeed components of the state ($x$) and their indistinguishable values are the corresponding components of $x'(t,~\tau)$.
\end{rem}

\begin{rem}\label{Remark2And3Include1}
If the state is observable, the orthogonal distribution $\OBS^\bot$ only consists of the null vector. Hence, the condition in (\ref{EquationConditionGeneral}) is trivially satisfied ($\xi^\alpha_i$ in (\ref{EquationConditionCoeffXi}) vanishes $\forall\alpha$ and $\forall i$). As a result, the first $m$ unknown inputs are locally identifiable. On the other hand, Theorem \ref{TheoremIdentifiabilityCan} states that the last $m_w-m$ unknown inputs are always unidentifiable. Therefore, from Theorem \ref{TheoremIdentifiabilityGeneral} and \ref{TheoremIdentifiabilityCan} we conclude that, if the state is observable, all the unknown inputs are locally identifiable if and only if $m=m_w$. This is precisely the statement of
Theorem \ref{TheoremIdentifiabilityObservable}.
\end{rem}

\begin{rem}\label{RemarkUnobservability}
State unobservability does not mean that the first $m$ unknown inputs are not locally identifiable. The condition in (\ref{EquationConditionGeneral}) could be honoured also in such a case (e.g., in the special case when all the quantities in (\ref{EquationConditionCoeffXi}) vanish for any $\xi\in\OBS^\bot$, i.e., when, for any $i=1,\ldots,m$, we have: $d\dt{\Li}_{g^0}\widetilde{h}_i\in\OBS$, and $d\Li_{g^k}\widetilde{h}_i\in\OBS$, for any $k=1,\ldots,m$.
\end{rem}



\subsection{Identifiability properties of $\E$ and of the original system}\label{SubSectionConditionSigma}

The analytical condition provided by Theorems \ref{TheoremIdentifiabilityObservable}, \ref{TheoremIdentifiabilityGeneral}, and \ref{TheoremIdentifiabilityCan} regards the system $\E$, which may differ from the original system because some of their unknown inputs, together with their time derivatives up to a given order, were included in the state ($\E$ may also differ from the original system because its unknown inputs could have been re-ordered). 
Let us suppose that this is the case and that the $j^{th}$ unknown input, together with its first $\mathcal{J}-1$ time derivatives, were included in the state of $\E$.

The condition of the above theorems actually provides the local identifiability of the $\mathcal{J}^{th}$ order time derivative of $w_j$ (i.e., $w^{(\mathcal{J})}_j:=\frac{d^{\mathcal{J}}w_j}{dt^{\mathcal{J}}}$). 
On the other hand,
the identifiability of the original unknown input and its first $\mathcal{J}-1$ time derivatives is trivially checked by verifying if their differentials belong to the observability codistribution $\OBS$, which is provided by Algorithm 9 in \cite{IF22}. In other words, it suffices to verify the condition ($l=0,1,\ldots,\mathcal{J}-1$.):

\begin{equation}\label{EquationConditionOriginalUI}
dw^{(l)}_j\in\OBS.
\end{equation}

\subsection{Identifiability of constant parameters}\label{SubSectionConditionConstant}

The constant parameters $Q_1,~\ldots,~Q_q$ belong to the state of $\E$. The condition that characterizes their local identifiability is trivially ($i=1,\ldots,q$):

\begin{equation}\label{EquationConditionConstant}
dQ_i\in\OBS.
\end{equation}

\section{HIV infection}\label{SectionHIV}

We investigate the identifiability properties of a simple ODE model widely used to describe
HIV dynamics in HIV-infected patients with antiretroviral treatment
\cite{Miao11,Villa19b,Per99,Per02,Chen08,Xia03,Moog10}.
The model is characterized by the following three equations and the following two outputs:

\begin{equation}\label{EquationHIVSystem}
\left\{\begin{array}{ll}
\dot{T}_U &= \lambda -\rho T_U -\eta(t) T_U V\\
\dot{T}_I &= \eta(t) T_U V -\delta T_I\\
\dot{V} &= N\delta T_I - cV\\
y &= [V,~T_U+T_I], \\
\end{array}\right.
\end{equation}
where, $T_U$ is the concentration of uninfected cells, $T_I$ the concentration of infected cells, and $V$ the viral load. The model is also characterized by the time varying parameter $\eta(t)$ and the five constant parameters $\lambda,~\rho,~\delta,~N,~c$. They are defined as follows:

\begin{itemize}

\item $\eta(t)$ is the infection rate, which is a function of the antiviral treatment efficacy. 
\item $\lambda$ is the source rate of uninfected cells. 
\item $\rho$ is the death rate of uninfected cells. 
\item $\delta$ is the death rate of infected cells. 
\item $N$ is the average number of virions produced by a single infected cell during its lifetime. 
\item $c$ is the clearance rate of free virions.

\end{itemize}

All these parameters are assumed to be unknown.

\vskip.2cm

To proceed, we need, first of all, to introduce a state that includes both the time varying quantities $T_U,~T_I,~V$, and the constant parameters (i.e., $\lambda,~\rho,~\delta,~N,~c$). We set:

\begin{equation}\label{EquationHIVState}
x=[T_U,~T_I,~V,~\lambda,~\rho,~\delta,~N,~c]^T
\end{equation}

The system is directly a special case of (\ref{EquationSystemDefinitionUIO}). In particular, 
$m_u=0$, $m_w=1$, $p=2$, $h_1(x)=V$, $h_2(x)=T_U+T_I$,

\begin{equation}\label{EquationHIVg0g1}
g^0=
\left[
\begin{array}{c}
 \lambda -\rho T_U\\
 -\delta T_I\\
N\delta T_I - cV\\
0\\
0\\
0\\
0\\
0\\
\end{array}
\right],~~
g^1=
\left[
\begin{array}{c}
 -T_UV\\
 T_UV\\
0\\
0\\
0\\
0\\
0\\
0\\
\end{array}
\right]
\end{equation}

\subsection{Observability analysis}\label{SubSectionHIVObs}

Algorithm 9 in \cite{IF22} provides, for this specific case, the following outcomes\footnote{A detailed derivation is available in Section 6.1 of \cite{arXivODE} (note that this derivation uses an equivalent version of Algorithm 9 in \cite{IF22}, which is Algorithm 1 in \cite{arXivODE}).}:

\begin{enumerate}

\item The system, denoted by $\E$ in this paper, coincides with the original system.

\item The unknown input degree of reconstructability of $\E$ is $m=m_w=1$ (system canonic with respect to its UI). In addition:

\begin{equation}\label{EquationHIVhtilde}
\widetilde{h}_1=\lambda-\rho T_U-\delta T_I.
\end{equation}

\item The observability codistribution

\[
\OBS=
\textnormal{span}\left\{
d h_1,~d h_2,~d h_3,~d \widetilde{h}_1,~dh_4,~dh_5,~dh_6
\right\},
\]
with: $h_3=N\delta T_I - cV$,
$h_4 =c (V c - N T_I \delta ) + \frac{N \delta  \rho   (T_I \delta  - \lambda  + T_U \rho  )}{\delta  - \rho}$, $h_5 =\frac{N \delta}{\rho-\delta}$, and $h_6 = c\frac{V c^2 \delta  - V c^2 \rho   - N T_I c \delta^2 + N T_I c \delta  \rho   + N T_I \delta^2 \rho   + N T_U \delta  \rho^2 - N \lambda  \delta  \rho}{\rho-\delta}$.

%
%
%
%
%
%
%
%
%


\end{enumerate}

The rank of $\OBS$ is $7$, which is smaller than the dimension of the state in (\ref{EquationHIVState}). Hence, the state is not observable. In particular, its orthogonal distribution is not empty. We have:

\begin{equation}\label{EquationHIVSymmetry}
\OBS^\bot=\textnormal{span}\left\{
\left[\begin{array}{l}
T_I\delta\\
-T_I\delta\\
0\\
0\\
0\\
\delta(\delta-\rho)\\
N\rho\\
0\\
\end{array}\right]\right\}.
\end{equation}

\subsection{Identifiability of the constant parameters}\label{SubSectionHIVIdentifiabilityConstant}

It is immediate to verify the following:

\[
d\lambda\in\OBS,~~
d\rho\in\OBS,~~
d\delta\notin\OBS,~~
dN\notin\OBS,~~
dc\in\OBS.
\]
(it suffices to check whether the above differentials are or not orthogonal to $\OBS^\bot$).
Therefore, regarding the constant parameters, two of them ($\delta$ and $N$) are unidentifiable, even locally.
This result contradicts the result available in the state of the art (e.g., see Section 6.2 of \cite{Miao11}).
~In Section \ref{SubSectionHIVComparisonSOTA}, we explicitly prove the validity of our result.

\subsection{Identifiability of the time varying parameter}\label{SubSectionHIVIdentifiabilityTV}

As $m=m_w$, all the identifiability properties are provided by Theorem \ref{TheoremIdentifiabilityGeneral}.
In accordance with Equation (\ref{EquationConditionGeneral}), we need to compute 
$\xi^0_1$, $\xi^1_1$, and $\nu^1_1$.
From (\ref{EquationHIVg0g1}) and (\ref{EquationHIVhtilde}) we obtain:
$\Li_{g^0}\widetilde{h}_1=T_I\delta^2 - \rho(\lambda - T_U\rho)$,
$\Li_{g^1}\widetilde{h}_1= -T_UV(\delta - \rho)$.
From Equation (\ref{EquationConditionCoeffXi}), by using the generator of $\OBS^\bot$ in (\ref{EquationHIVSymmetry}) for $\xi$, we obtain:

\begin{equation}\label{EquationHIVCoeffXi}
\xi^0_1=T_I\delta(\delta - \rho)^2,~~
\xi^1_1=-V\delta(T_I + T_U)(\delta - \rho).
\end{equation}

Let us compute $\nu^1_1$. From (\ref{EquationConditionMu}), (\ref{EquationHIVg0g1}), and (\ref{EquationHIVhtilde}) we obtain $\mu^1_1=-T_UV(\delta - \rho)$. Hence, for its inverse:

\begin{equation}\label{EquationHIVNu}
\nu^1_1=-\frac{1}{T_UV(\delta - \rho)}.
\end{equation}

The condition in (\ref{EquationConditionGeneral}) is not honoured (for both $\alpha=0,1$).
This means that the time varying parameter $\eta(t)$ is not identifiable, even locally.
This result contradicts the result available in the state of the art
(e.g., see Section 6.2 of \cite{Miao11} ).
~In Section \ref{SubSectionHIVComparisonSOTA}, we explicitly prove the validity of our result.

We apply the second part of Theorem \ref{TheoremIdentifiabilityGeneral} in order to determine infinitely many values of the non identifiable parameters which are indistinguishable from the true values. This regards both the constant and the time varying parameters, according to Remark \ref{RemarkContinuousConstant}.
To determine these values,
we need to solve the system of differential equations in (\ref{EquationConditionDiffEqSystem}), for the specific case.
By substituting the expressions in (\ref{EquationHIVCoeffXi}) and (\ref{EquationHIVNu}) in (\ref{EquationConditionChij}) we obtain:

\[
\chi_1=\frac{T_I\delta(\delta - \rho)}{T_UV}
-\frac{\delta(T_I + T_U)}{T_U}\eta
\]

As a result, the system of differential equations in (\ref{EquationConditionDiffEqSystem}) becomes:

\begin{equation}\label{EquationHIVDiffEqSystem}
\left\{\begin{array}{ll}
  \frac{dT_U'}{d\tau} &=  T_I'\delta'\\
  \frac{dT_I'}{d\tau} &=  -T_I'\delta'\\
  \frac{dV'}{d\tau} &=  0\\
  \frac{d\lambda'}{d\tau} &=  0\\
  \frac{d\rho'}{d\tau} &=  0\\
  \frac{d\delta'}{d\tau} &=  \delta'(\delta'-\rho')\\  
  \frac{dN'}{d\tau} &=   N'\rho'\\
  \frac{dc'}{d\tau} &=  0\\
  \frac{d\eta'}{d\tau} &=  \frac{T_I'\delta'(\delta' - \rho')}{T_U'V'}
-\frac{\delta'(T_I' + T_U')}{T_U'}\eta'\\
 T_U'(t,~0)&=T_U(t), ~T_I'(t,~0)= T_I(t),\\
 V'(t,~0)&= V(t),\\
  \lambda'(0)= \lambda, ~\rho'(0)&= \rho, ~\delta'(0)= \delta, ~N'(0)= N,\\ ~c'(0)= c,~ \eta'(t,~0)&= \eta(t) \\
\end{array}\right.
\end{equation}
These equations can be solved analytically. We obtain for $x'(t,~\tau)$:

\begin{equation}\label{EquationHIVIndistinguishableStates}
\left\{\begin{array}{ll}
T_U'(t,~\tau)&=T_U(t)+ T_I(t)-\frac{ T_I(t)}{ \rho}\left(
 \delta e^{- \rho\tau}+ \rho- \delta
\right)\\
 T_I'(t,~\tau)&=\frac{ T_I(t)}{ \rho}\left(
 \delta e^{- \rho\tau}+ \rho- \delta
\right)\\
V'(t,~\tau)&=V(t)\\
\lambda'(\tau)&= \lambda\\
\rho'(\tau)&= \rho\\
\delta'(\tau)&= \delta \rho\left/
\left(( \rho- \delta)e^{ \rho\tau}+ \delta\right)\right.\\
N'(\tau)&= Ne^{ \rho\tau}\\
c'(\tau)&= c \\
\end{array}\right.
\end{equation}

Regarding the unknown input, we obtain:

\begin{equation}\label{EquationHIVIndistinguishableUI}
\eta'(t,~\tau)=
\end{equation}
{\scriptsize
\[
\frac{ \eta(t) T_U(t)  V(t)  \rho e^{ \rho \tau} +\left[T_I(t)  \delta^2-  T_I(t)  \delta  \rho-  \eta(t) T_U(t)  V(t)  \delta  \right] \left[e^{ \rho \tau}-1 \right]}{ V(t) \left[ T_I(t)\delta +  T_U(t)  \rho \right] e^{ \rho \tau} -  V(t)T_I(t)  \delta}
\]
}

We obtained the following fundamental result. 
Let us denote by $x(t)=[T_U(t),~T_I(t),~V(t),~\lambda,~\rho,~\delta,~N,~c]^T$ the true state at the time $t$, and by $\eta(t)$ the true unknown input at the time $t$. Then, $x'(t,~\tau)$ given in (\ref{EquationHIVIndistinguishableStates}) and $\eta'(t,~\tau)$ given in (\ref{EquationHIVIndistinguishableUI}) are indistinguishable from $x(t)$ and $\eta(t)$, for any $\tau\in\R$.


\subsection{Comparison with the state of the art}\label{SubSectionHIVComparisonSOTA}

As we aforementioned, our results contradict the result available in the state of the art for exactly the same ODE model characterized by 
(\ref{EquationHIVSystem}) (e.g., see Section 6.2 of \cite{Miao11}).
~Specifically, in contrast with the state of the art, we obtained that both the state and the parameters, cannot be recovered from the outputs. In addition, we also obtained a one dimensional set of indistinguishable states, $x'(t,~\tau)$, and a one parameter set of functions $\eta'(t,~\tau)$ such that, by only using the outputs, it is not possible to distinguish among the elements of these sets.
Because of this discrepancy, we wish to provide here an easy and direct proof that  
unequivocally shows the unidentifiability of the system, in contrast with the results available in the state of the art. In particular, 
the validity of our results can be checked by a trivial substitution. The interested reader can also find some further details in Section 6.6 of \cite{arXivODE} and in \cite{arXivErratum}.

First, from (\ref{EquationHIVIndistinguishableStates}), it is trivial to check that the outputs $y_1(t)=T_U(t)+T_I(t)$ and $y_2(t)=V(t)$ coincide with $T_U'(t,~\tau)+T_I'(t,~\tau)$ and $V'(t,~\tau)$, and this holds for any $\tau\in\R$.

Second, we compute the time derivatives of $T_U'(t,~\tau)$, $T_I'(t,~\tau)$, and $V'(t,~\tau)$. We use their expression in (\ref{EquationHIVIndistinguishableStates}) and we use 
(\ref{EquationHIVSystem}). By also using the expression of $\eta'(t,~\tau)$ in (\ref{EquationHIVIndistinguishableUI}) we finally obtain:

\begin{equation}\label{EquationHIVDynamicsTransformed}
\left\{\begin{array}{ll}
\frac{d}{dt}T_U' &=\lambda'-\rho'T_U'-\eta'T_U'V'\\
\frac{d}{dt}T_I' &=\eta'T_U'V'-\delta'T_I'\\
\frac{d}{dt}V' &=N'\delta'T_I'-c'V'\\
\end{array}\right.
\end{equation}

In other words, the new three states, satisfy exactly the same dynamical equations in (\ref{EquationHIVSystem}) with the new parameters.

We conclude that
$T_U'(t, ~\tau)$, $T_I'(t, ~\tau)$, $V'(t, ~\tau)$, 
$\eta'(t, ~\tau)$, $\lambda'(\tau)$, $\delta'(\tau)$, $\rho'(\tau)$, $N'(\tau)$, and $c'(\tau)$ cannot be distinguished from $T_U(t)$, $T_I(t)$, $V(t)$, $\eta(t)$, $\lambda$, $\delta$, $\rho$, $N$, and $c$, respectively. 
Regarding the model parameters, 
as $\eta'(t,~\tau)\neq\eta(t)$, $\delta'(\tau)\neq\delta$, and $N'(\tau)\neq N$, we conclude that $\eta(t)$, $\delta$, and $N$ are not identifiable. Finally, we remark that $\tau$ can take infinitesimal values and that, in the limit of $\tau\rightarrow0$, we have: 
$T_U'(t,~\tau)\rightarrow T_U(t)$, $T_I'(t,~\tau)\rightarrow T_I(t)$, $V'(t,~\tau)\rightarrow V(t)$,
$\eta'(t,~\tau)\rightarrow\eta(t)$, $\delta'(\tau)\rightarrow\delta$, and $N'(\tau)\rightarrow N$. Therefore, we conclude that the three parameters $\eta(t)$, $\delta$, and $N$ are not even locally identifiable.

We wish to emphasize that our claim was here easily proved by a simple substitution. What is not trivial is the analytic derivation of the expressions in 
(\ref{EquationHIVIndistinguishableStates}), and (\ref{EquationHIVIndistinguishableUI}). These expressions were analytically determined by using the very powerful and general result stated by Theorem \ref{TheoremIdentifiabilityGeneral}. Thanks to this result, the derivation of the expressions in 
(\ref{EquationHIVIndistinguishableStates}), and (\ref{EquationHIVIndistinguishableUI}) becomes trivial in the sense that it follows the steps of an automatic procedure that requires no inventiveness from the user.

\subsection{Numerical Results}\label{SubSectionHIVNumerical}

In this section, we explicitly provide the indistinguishable sets  obtained in Section \ref{SubSectionHIVIdentifiabilityTV}. We refer
to the same data available in the literature (e.g., see the electronic supplementary material of \cite{Villa19b}).
~The data set is characterized as follows:

\begin{itemize}

\item The time interval is $\mathcal{I}=[0,~201]s$.

\item $T_U(0)=600$, $T_I(0)=0$, $V(0)=10^5$, $\lambda=36$, $\rho=0.108$, $\delta=0.5$, $N=10^3$, and $c=3$.

\item The time varying parameter is:
\begin{equation}\label{EquationHIVEtaDataSet}
\eta(t)=k\left(
1-0.9\cos\left(\frac{\pi}{1000}t\right)
\right),
\end{equation}
with $k=0.00009$.

\end{itemize}

\begin{figure}[htbp]
\begin{center}
\includegraphics[width=.8\columnwidth]{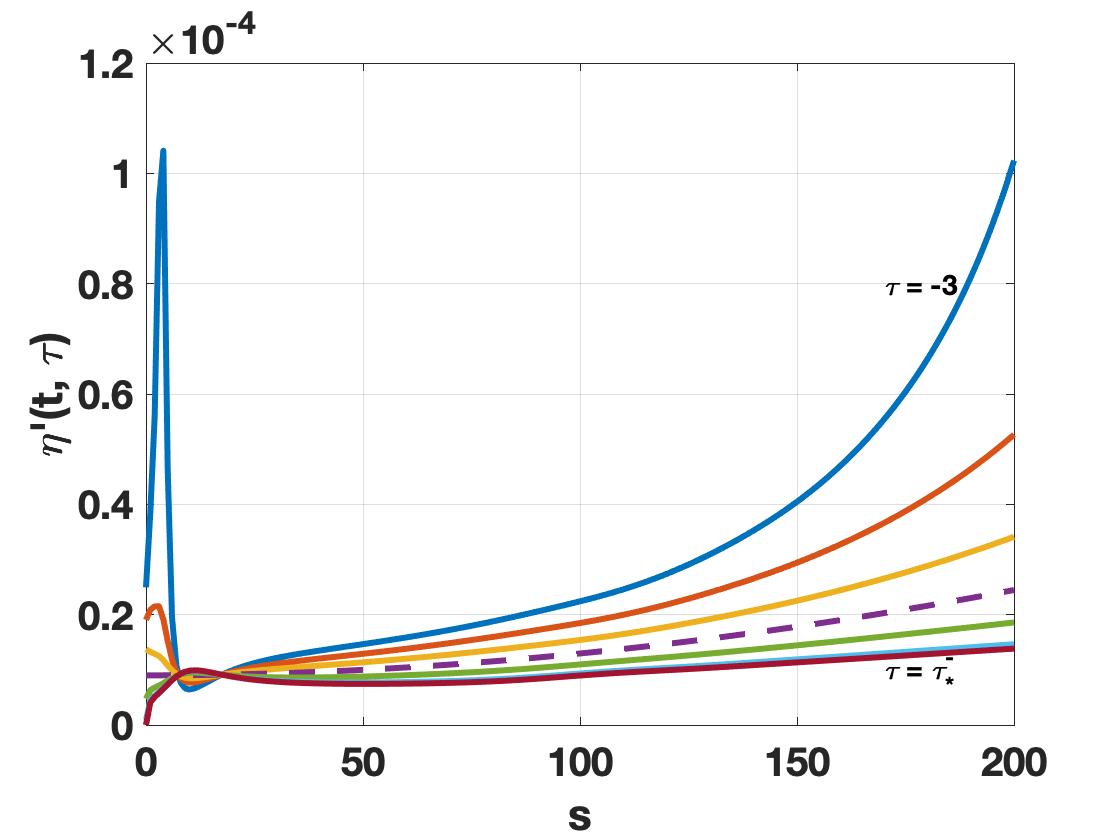}
\caption{Several indistinguishable profiles for the time varying parameter $\eta'(t,~\tau)$ for several values of the parameter $\tau$ ranging from $-3$ up to $\tau_*\cong2.2532$. The dashed purple line is the profile given by (\ref{EquationHIVEtaDataSet}), i.e., $\eta'(t,~\tau)$ for $\tau=0$.} \label{FigHIVEta}
\end{center}
\end{figure}

\begin{figure}[htbp]
\begin{center}
\begin{tabular}{cc}
\includegraphics[width=.45\columnwidth]{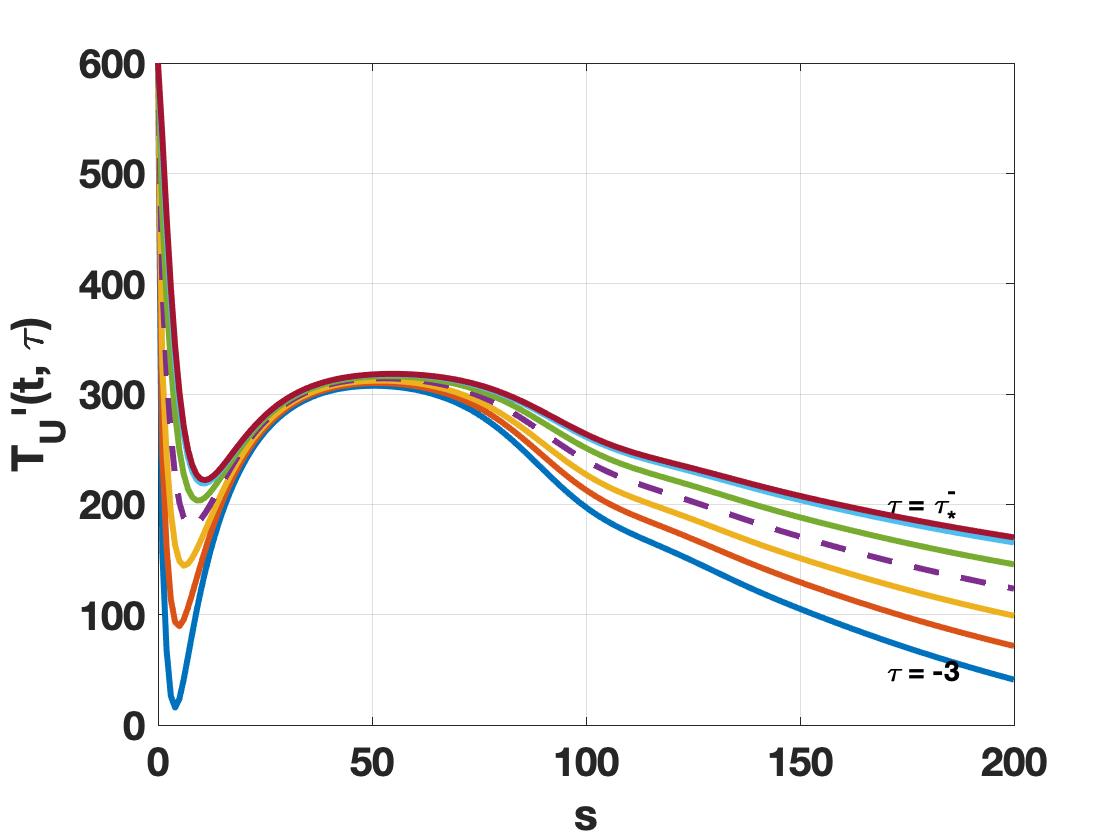}&\includegraphics[width=.45\columnwidth]{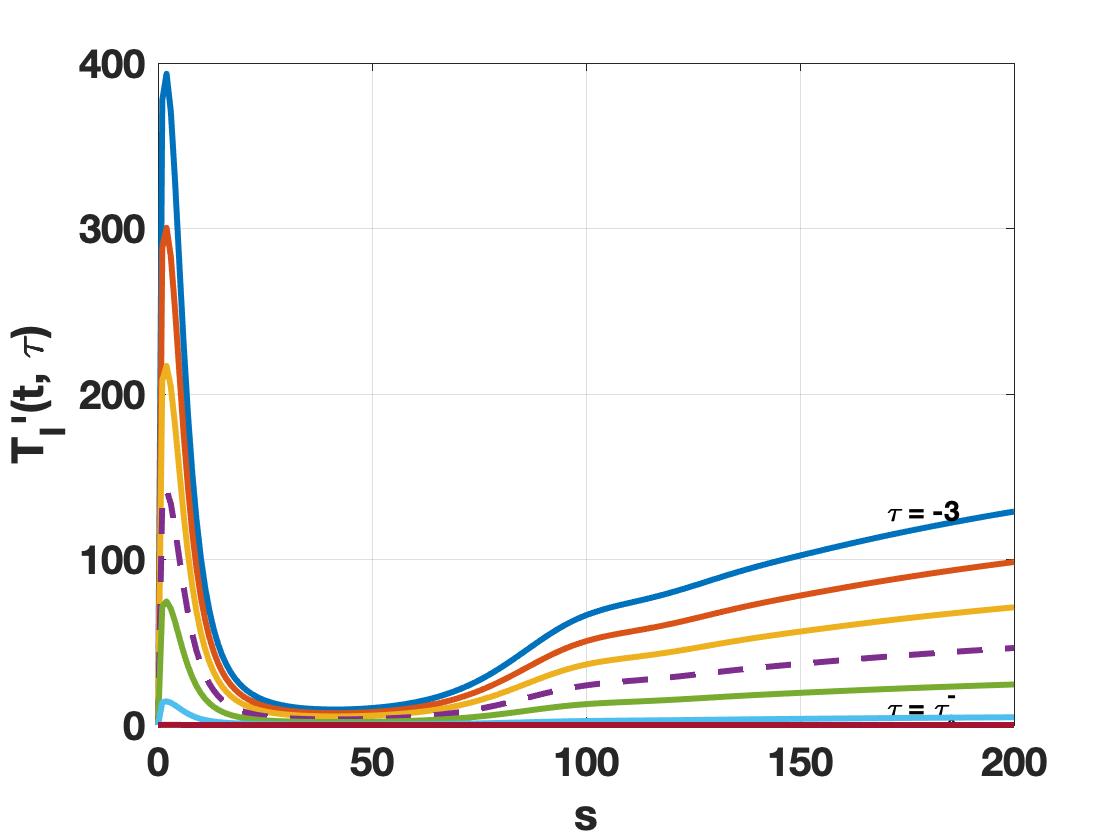}\\
\end{tabular}
\caption{The profiles $T_U'(t,~\tau)$ (left side) and $T_I'(t,~\tau)$ (right side) for several values of the parameter $\tau$ ranging from $-3$ up to $\tau_*\cong2.2532$. The dashed purple line is the profile for $\tau=0$.} \label{FigHIVTuTi}
\end{center}
\end{figure}

%

%
%
%
%

We provide indistinguishable states and indistinguishable unknown inputs for several values of the parameter $\tau\in\mathcal{U}$, where $\mathcal{U}$ is an interval of $\R$ that includes $\tau=0$ and such that, for any $\tau\in\mathcal{U}$, the quantities in (\ref{EquationHIVIndistinguishableStates}) and (\ref{EquationHIVIndistinguishableUI}) have physical meaning, as explained below.

From the expression of $\delta'(\tau)$ in (\ref{EquationHIVIndistinguishableStates}) and the values of our data set, we obtain that

\[
\lim_{\tau\rightarrow\tau_*^-}\delta'(\tau)=+\infty,~~
\lim_{\tau\rightarrow\tau_*^+}\delta'(\tau)=-\infty
\]
with $\tau_*:=\frac{\log\left(\frac{\delta}{\delta-\rho}\right)}{\rho}\cong 2.2532$.
In addition, we found that, for $\tau=-4$ the function $\eta'(t,~\tau)$ in (\ref{EquationHIVIndistinguishableUI}) takes negative values for $t\in[2,~6]s$.
Finally, we verified that in the interval
\[
\mathcal{U}:=\left[-3,~\frac{}{}\tau_*\right)
\]
all the quantities $\delta'(\tau)$, $N'(\tau)$, $\eta'(t,~\tau)$, $T_U'(t,~\tau)$, and $T_I'(t,~\tau)$ take positive values (the last three on the entire time interval $\mathcal{I}$).
Hence, we are allowed to use any $\tau\in\mathcal{U}$.


Figure \ref{FigHIVEta} displays the profiles of $\eta'(t,~\tau)$. The profile for $\tau=0$ (purple dashed line) is precisely the one set in (\ref{EquationHIVEtaDataSet}). By varying $\tau$ we obtain a significant change of the profile meaning that this parameter is "strongly" unidentifiable.
In addition, by varying $\tau$, the two constant parameters that are not identifiable (i.e., $\delta$ and $N$) significantly change. For instance, in accordance with (\ref{EquationHIVIndistinguishableStates}),
for $\tau=-3$ and $\tau=\tau_*$ we obtain the following changes:
\[
\begin{array}{lll}
\delta=\delta'(0)=0.5 & \hskip-.2cm\rightarrow \delta'(-3)\cong0.25 &\hskip-.2cm\rightarrow \delta'(\tau\rightarrow\tau_*^-)=+\infty\\
N=N'(0)=10^3& \hskip-.2cm\rightarrow N'(-3)\cong 723 &\hskip-.2cm\rightarrow N'(\tau_*)\cong 1276.\\
\end{array}
\]
The variation of $\delta'$ is even divergent.

Figure \ref{FigHIVTuTi} displays the profiles of $T_U'(t,~\tau)$ and $T_I'(t,~\tau)$. They also change significantly by varying $\tau$ meaning that also $T_U'(t,~\tau)$ and $T_I'(t,~\tau)$ are "strongly" unobservable.

\subsection{Minimal external information}\label{SectionHIVMinimalExtInfo}

Our analysis on the model characterized by (\ref{EquationHIVSystem}) showed that only three parameters are identifiable. They are $\lambda,~\rho$, and $c$. The remaining constant parameters, $\delta$, and $N$, and the time varying parameter, $\eta(t)$, cannot be identified.
We conclude this study by answering to the following fundamental practical question: {\it which measurements must be added to the model to make all the model parameters locally identifiable?}

It is possible to prove that, by also including the output $y_3=T_U$ (or $y_3=T_I$), the state becomes observable and all the parameters identifiable. This can be proved by repeating the analysis for the new system characterized by this further output.
On the other hand, from our analysis, it is immediate to obtain a more interesting result.
It actually suffices to have this further output at a single time $t^*\in\mathcal{I}$ (and not necessarily on the entire time interval) to make possible the identification of all the parameters (and the observability of the state). Indeed, by knowing $T_U$ (or equivalently $T_I=y_2-T_U$) at a single time, we can easily determine the value of $\tau$ from the expression of $T_U'$ (or $T_I'$) in (\ref{EquationHIVIndistinguishableStates}).

%

\section{Compartmental SEIAR model of the epidemic dynamics of Covid-19}\label{SectionCovid}

We consider a simple model that generalizes the SEIR model commonly used for virus disease. 
This more general ODE model includes a further compartment, $A(t)$, that is the asymptomatic infected, \cite{SEIARPribylova}. 
The other compartments are the same of a SEIR model, i.e.: S, Susceptible; E, Exposed; I, Infected (symptomatic infected); R, Removed (i.e., healed or dead).
The differential equations of this model 
describe the movement of individuals of the
population between the five aforementioned classes, i.e.: $S(t),~E(t),~I(t),~A(t)$ and $R(t)$. 
The dynamics and the outputs are given by the following set of equations (e.g., see \cite{SEIARPribylova}):

\begin{equation}\label{EquationCovidSystem}
\left\{\begin{array}{ll}
\dot{S} &= -\beta S(I+A) \\
\dot{E} &= \beta S(I+A) -\gamma E\\
\dot{I} &= \gamma pE-\mu_1I\\
\dot{A} &= \gamma(1-p)E-\mu_2A\\
\dot{R} &= \mu_1I+\mu_2A\\
y &= [I, ~A, ~S+E+R], \\
\end{array}\right.
\end{equation}

This model includes the following constant parameters, which are assumed to be unknown:

\begin{itemize}

\item $\mu_1$ and $\mu_2$ that are the remove rates from the infected symptomatic individuals ($I(t)$) and the infected asymptomatic individuals ($A(t)$), respectively.

\item $\gamma$ that is the rate
at which the exposed individuals ($E(t)$) are infected. 

\item $p$ that is the probability that infected individuals are symptomatic (and $1-p$  is the probability that infected individuals are asymptomatic).

\end{itemize}

The model also includes the following unknown time varying parameter:

\begin{itemize}

\item $\beta=\beta(t)$ that is the probability of disease transmission in a single contact
times the average number of contacts per person, due to contacts between an individual of the class $S$ and an individual that belongs to one of the two classes $I$ and $A$. This parameter is assumed to be time varying as it can be modified by government measures (e.g., , using masks, closing schools, remote working).

\end{itemize}

We assume that we can measure the infected symptomatic individuals ($I(t)$), and the infected asymptomatic individuals ($A(t)$).
Hence, our model is characterized by the first two outputs in (\ref{EquationCovidSystem}). 
Finally, as we know the total population, $S(t)+E(t)+I(t)+A(t)+R(t)$, by removing from this the above two measurements, we obtain the third output in (\ref{EquationCovidSystem}), i.e., $S(t)+E(t)+R(t)$.

\vskip.2cm

To proceed, we need, first of all, to introduce a state that includes both the time varying quantities (i.e., $S,~E,~I,~A,~R$) and the constant parameters (i.e., $\mu_1, ~\mu_2, ~\gamma,~p$). We set:

\begin{equation}\label{EquationCovidState}
x=[S,~E,~I, ~A, ~R, ~\mu_1, ~\mu_2, ~\gamma,~p]^T
\end{equation}

The system is directly a special case of (\ref{EquationSystemDefinitionUIO}). In particular, 
$m_u=0$, $m_w=1$, $p=3$, $h_1(x)=I$, $h_2(x)=A$, $h_3(x)=S+E+R$,

\begin{equation}\label{EquationCovidg0g1}
g^0=
\left[
\begin{array}{c}
0\\
-\gamma E\\
\gamma pE-\mu_1I\\
\gamma(1-p)E-\mu_2A\\
\mu_1I+\mu_2A\\
0\\
0\\
0\\
0\\
\end{array}
\right],~~
g^1=
\left[
\begin{array}{c}
 -S(I+A)\\
 S(I+A)\\
0\\
0\\
0\\
0\\
0\\
0\\
0\\
\end{array}
\right]
\end{equation}

\subsection{Observability analysis}\label{SubSectionCovidObs}

Algorithm 9 in \cite{IF22} provides, for this specific case, the following outcomes\footnote{A detailed derivation is available in Section 7.1 of \cite{arXivODE} (note that this derivation uses an equivalent version of Algorithm 9 in \cite{IF22}, which is Algorithm 1 in \cite{arXivODE}).}:

\begin{enumerate}

\item The system, denoted by $\E$ in this paper, coincides with the original system.

\item The unknown input degree of reconstructability of $\E$ is $m=m_w=1$ (system canonic with respect to its UI). In addition:

\begin{equation}\label{EquationCovidhtilde}
\widetilde{h}_1=\gamma pE-\mu_1I.
\end{equation}

\item The observability codistribution\\
$
\OBS=
\textnormal{span}\left\{
d h_1,~d h_2,~d h_3,~dh_4,~dh_5,~dh_6,~d \widetilde{h}_1
\right\}$,
where $h_1(x)=I$, $h_2(x)=A$, $h_3(x)=S+E+R$, $h_4(x)=-A\mu_2-E\gamma(p-1)$,
$h_5(x)=\frac{1-p}{p}$, and $h_6(x)=\mu_2(A\mu_2 + E\gamma(p - 1)) + \mu_1\frac{1-p}{p}(E\gamma p-I\mu_1)$.
\end{enumerate}


The rank of $\OBS$ is $7$, which is smaller than the dimension of the state in (\ref{EquationCovidState}). Hence, the state is not observable. In particular, its orthogonal distribution is not empty. We have:

\begin{equation}\label{EquationCovidSymmetry}
\OBS^\bot=\textnormal{span}\left\{
\left[\begin{array}{c}
1\\
0\\
0\\
0\\
-1\\
0\\
0\\
0\\
0\\
\end{array}\right],~
\left[\begin{array}{c}
E\\
-E\\
0\\
0\\
0\\
0\\
0\\
\gamma\\
0\\
\end{array}\right]
\right\}.
\end{equation}

\subsection{Identifiability of the constant parameters}\label{SubSectionCovidIdentifiabilityConstant}

It is immediate to verify the following:

\[
d\mu_1\in\OBS,~~
d\mu_2\in\OBS,~~
d\gamma\notin\OBS,~~
dp\in\OBS.
\]
(it suffices to check whether the above differentials are or not orthogonal to $\OBS^\bot$).
Therefore, regarding the constant parameters, one of them ($\gamma$) is unidentifiable, even locally.

\subsection{Identifiability of the time varying parameter}\label{SubSectionCovidIdentifiabilityTV}

As $m=m_w$, all the identifiability properties are provided by Theorem \ref{TheoremIdentifiabilityGeneral}.
In accordance with Equation (\ref{EquationConditionGeneral}), we need to compute $\nu^1_1$, and, for both the generators of $\OBS^\bot$ in (\ref{EquationCovidSymmetry}), the coefficients 
$\xi^0_1$, $\xi^1_1$.

Let us compute $\nu^1_1$. From (\ref{EquationConditionMu}), (\ref{EquationCovidg0g1}), and (\ref{EquationCovidhtilde}) we obtain $\mu^1_1=S\gamma p(A + I)$. Hence, for its inverse:

\begin{equation}\label{EquationCovidNu}
\nu^1_1=\frac{1}{S\gamma p(A + I)}.
\end{equation}

From (\ref{EquationCovidg0g1}) and (\ref{EquationCovidhtilde}) we obtain:

\[
\Li_{g^0}\widetilde{h}_1=\mu_1(I\mu_1 - E\gamma p) - E\gamma^2p,~~
\Li_{g^1}\widetilde{h}_1= S\gamma p(A + I).
\]

From Equation (\ref{EquationConditionCoeffXi}), by using the first generator of $\OBS^\bot$ in (\ref{EquationCovidSymmetry}) for $\xi$, we obtain:

\begin{equation}\label{EquationCovidCoeffXi}
\xi^0_1=0,~~
\xi^1_1=\gamma p(A + I).
\end{equation}
The condition in (\ref{EquationConditionGeneral}) is not honoured for $\alpha=1$. It is unnecessary to repeat the computation with the second generator of $\OBS^\bot$ in (\ref{EquationCovidSymmetry}). We conclude that
the time varying parameter $\beta(t)$ is not identifiable, even locally.

We apply the second part of Theorem \ref{TheoremIdentifiabilityGeneral} in order to determine values of the non identifiable parameters which are indistinguishable from the true values. This regards both the constant and the time varying parameter, according to Remark \ref{RemarkContinuousConstant}. As in this case we have two independent generators of $\OBS^\bot$, we can determine two independent sets.

\subsubsection{First set} Let us start with the first generator.
We provide the corresponding system of differential equations in (\ref{EquationConditionDiffEqSystem}).
By using (\ref{EquationCovidCoeffXi}) and (\ref{EquationCovidNu}) in (\ref{EquationConditionChij}) we obtain:

\[
\chi_1 = -~\nu^1_1 ( \xi^0_1 + \xi^1_1 w_1 )=
-\frac{1}{S} \beta
\]

As a result, the system of differential equations in (\ref{EquationConditionDiffEqSystem}) becomes:

\begin{equation}\label{EquationCovidDiffEqSystem1}
\left\{\begin{array}{ll}
  \frac{dS'}{d\tau} &= 1 \\
  \frac{dE'}{d\tau} &= 0 \\
  \frac{dI'}{d\tau} &=  0\\
  \frac{dA'}{d\tau} &=  0\\
  \frac{dR'}{d\tau} &=  -1\\
  \frac{d\mu_1'}{d\tau} &= 0 \\  
  \frac{d\mu_2'}{d\tau} &=  0\\  
  \frac{d\gamma'}{d\tau} &= 0 \\
  \frac{dp'}{d\tau} &=  0\\
   \frac{d\beta'}{d\tau} &=  -\frac{\beta'}{S'}\\
S'(t,~0)&=S(t), ~E'(t,~0)= E(t), ~I'(t,~0)= I(t), \\
A'(t,~0)&= A(t), ~R'(t,~0)= R(t), ~\mu_1'(0)= \mu_1,\\
 \mu_2'(0)&= \mu_2, ~\gamma'(0)= \gamma, ~p'(0)= p \\
 \beta'(t,~0)&= \beta(t) \\
\end{array}\right.
\end{equation}
These equations can be solved analytically. We obtain for $x'(t,~\tau)$:

\begin{equation}\label{EquationCovidIndistinguishableStates1}
\left\{\begin{array}{ll}
S'(t,~\tau)&=S(t)+\tau\\
E'(t,~\tau)&=E(t)\\
I'(t,~\tau)&=I(t)\\
A'(t,~\tau)&=A(t)\\
R'(t,~\tau)&=R(t)-\tau\\
~\mu_1'(\tau)&= \mu_1\\
\mu_2'(\tau)&= \mu_2\\
\gamma'(\tau)&= \gamma\\
p'(\tau)&= p \\
\end{array}\right.
\end{equation}


Regarding the unknown input, we obtain:

\begin{equation}\label{EquationCovidIndistinguishableUI1}
\beta'(t,~\tau)=\frac{\beta(t)S(t)}{S(t)+\tau}.
\end{equation}

\subsubsection{Second set} Let us use now the second generator.
We provide the corresponding system of differential equations in (\ref{EquationConditionDiffEqSystem}). First, instead of (\ref{EquationCovidCoeffXi}) we have:
\[
\xi^0_1=-\gamma^2~p~E,~~~
\xi^1_1=\gamma p (A +  I)(E+S).
\]

and, by substituting in (\ref{EquationConditionChij}), with $\nu^1_1$ given in (\ref{EquationCovidNu}), we obtain:

\[
\chi_1 = -~\nu^1_1 ( \xi^0_1 + \xi^1_1 w_1 )=
\frac{\gamma~E}{S (A +  I)}-\frac{E+S}{S}  \beta
\]

As a result, the system of differential equations in (\ref{EquationConditionDiffEqSystem}) becomes:

\begin{equation}\label{EquationCovidDiffEqSystem2}
\left\{\begin{array}{ll}
  \frac{dS'}{d\tau} &= E' \\
  \frac{dE'}{d\tau} &= -E' \\
  \frac{dI'}{d\tau} &=  0\\
  \frac{dA'}{d\tau} &=  0\\
  \frac{dR'}{d\tau} &=  0\\
  \frac{d\mu_1'}{d\tau} &= 0 \\  
  \frac{d\mu_2'}{d\tau} &=  0\\  
  \frac{d\gamma'}{d\tau} &= \gamma' \\
  \frac{dp'}{d\tau} &=  0\\
   \frac{d\beta'}{d\tau} &=  \frac{\gamma'~E'}{S' (A' +  I')}-\frac{(E'+S')}{S'}  \beta'\\
S'(t,~0)&=S(t), ~E'(t,~0)= E(t), ~I'(t,~0)= I(t), \\
A'(t,~0)&= A(t), ~R'(t,~0)= R(t), ~\mu_1'(0)= \mu_1,\\
 \mu_2'(0)&= \mu_2, ~\gamma'(0)= \gamma, ~p'(0)= p \\
 \beta'(t,~0)&= \beta(t) \\
\end{array}\right.
\end{equation}
These equations can be solved analytically. We obtain for $x'(t,~\tau)$:

\begin{equation}\label{EquationCovidIndistinguishableStates2}
\left\{\begin{array}{ll}
S'(t,~\tau)&=S(t)+E(t)~(1-e^{-\tau})\\
E'(t,~\tau)&=E(t)~e^{-\tau}\\
I'(t,~\tau)&=I(t)\\
A'(t,~\tau)&=A(t)\\
R'(t,~\tau)&=R(t)\\
~\mu_1'(\tau)&= \mu_1\\
\mu_2'(\tau)&= \mu_2\\
\gamma'(\tau)&= \gamma~e^{\tau}\\
p'(\tau)&= p \\
\end{array}\right.
\end{equation}

Regarding the unknown input, we obtain:

\begin{equation}\label{EquationCovidIndistinguishableUI2}
\beta'(t,~\tau)=\frac{\gamma E(t)(1-e^\tau)-S(t)\beta(t)\left[A(t)+I(t)\right]}{\left[A(t)+I(t)\right]\left[E(t)-(E(t)+S(t))e^\tau\right]}.
\end{equation}

\vskip.2cm
It is immediate to test the validity of our results by proceeding exactly as in Section \ref{SubSectionHIVComparisonSOTA}, for the two previous sets, i.e., by using the solution given in (\ref{EquationCovidIndistinguishableStates1}), (\ref{EquationCovidIndistinguishableUI1}) and the solution given in (\ref{EquationCovidIndistinguishableStates2}), (\ref{EquationCovidIndistinguishableUI2}). In Section 7.6 of \cite{arXivODE}, we also provide numerical results obtained with these solutions.

\section{Visual-inertial sensor fusion}\label{SectionVISFM}

We provide a further application in a completely different domain.
We refer to the visual inertial sensor fusion problem. 
In contrast with the viral systems studied in the previous sections, this system also includes the presence of known inputs ($m_u\neq0$).
For the sake of clarity, we restrict our analysis to a $2D$ environment. The $3D$ case differs only in a more laborious calculation.

\begin{figure}[htbp]
\begin{center}
\includegraphics[width=.6\columnwidth]{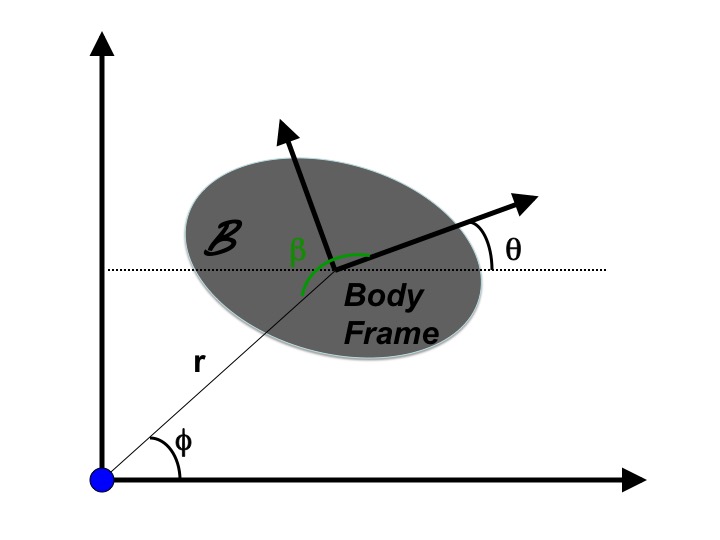}
\caption{The global frame, the body frame and the observation provided by the visual sensor (the angle $\beta$).} \label{FigVISFM}
\end{center}
\end{figure}

We consider a rigid body ($\mathcal{B}$) equipped with a visual sensor and a gyroscope. The body moves on a plane. The visual sensor provides the bearing angle of the features in its own local frame. The gyroscope provides the angular speed (which is a scalar, in $2D$).
We assume that the local frames of the two sensors coincide and we call this common frame, the body frame. In addition, we assume that the gyroscope measurements are unbiased.
Figure \ref{FigVISFM} depicts our system.


\noindent It is very convenient to work in polar coordinates. 
Hence, we define the state $[r,~\phi,~v,~\alpha,~\theta]^T$, where $r$ and $\phi$ characterize the body position, $\theta$ its orientation (see Fig. \ref{FigVISFM} for an illustration), and $v$ and $\alpha$ the body speed in polar coordinates. In particular,  $v= \sqrt{v_x^2+v_y^2}$ and $\alpha= \arctan\left(\frac{v_y}{v_x}\right)$, where $[v_x,~v_y]^T$ is the body speed in Cartesian coordinates. The dynamics are:

\begin{equation}\label{EquationVISFMDynamics}
\left[\begin{array}{ll}
  \dot{r} ~&= v \cos(\alpha-\phi)\\
  \dot{\phi} ~&= \frac{v}{r} \sin(\alpha-\phi)\\
  \dot{v} ~&= A_x \cos(\alpha-\theta)  +  A_y \sin(\alpha-\theta)\\
 \dot{\alpha}~&= -\frac{A_x}{v} \sin(\alpha-\theta)  + \frac{A_y}{v} \cos(\alpha-\theta)\\
  \dot{\theta} ~~&= \omega \\
\end{array}\right.
\end{equation}
where $[A_x, ~A_y]^T$ is the body acceleration in the body frame and $\omega$ the angular speed.
For the sake of brevity, we assume that the visual sensor observes a single point feature and that this feature is positioned at the origin of the global frame. In other words, the visual sensor provides the angle $\beta=\pi-\theta+\phi$. Hence, we can perform the observability and the identifiability  analysis by using the output (we ignore $\pi$):

\begin{equation}\label{EquationVISFMOutput}
y=h(x)=\phi-\theta
\end{equation}

$A_x$, $A_y$ are unknown and $\omega$ is known.
By comparing (\ref{EquationVISFMDynamics}) and (\ref{EquationVISFMOutput}) with (\ref{EquationSystemDefinitionUIO}) we have: $n=5$, $m_u=1$, $m_w=2$,
$p=1$,
$u_1=\omega$, $w_1=A_x$, $w_2=A_y$.
As $p=1<2=m_w$ the system is certainly not in canonical form with respect to $w_1$ and $w_2$. We even do not know if it is canonic with respect to them.

\subsection{Observability analysis}\label{SubSectionVISFMStateObservability}

By applying Algorithm 9 in \cite{IF22} we obtain the system $\E$, its unknown input degree of reconstructability $m$, the functions $\widetilde{h}_1,\ldots,\widetilde{h}_m$, and the observability codistribution of $\E$ (denoted by $\OBS$)\footnote{A detailed derivation is available in Section 9.2.2 of \cite{IF22}. Note that this derivation is also available in the arXiv version of \cite{IF22} (where, however, instead of Algorithm 9 in \cite{arXivODE}, it is used an equivalent version of it).}. 
In particular, $\E$ differs from the original system because its state also includes the second unknown input ($A_y$). Hence, the state that defines $\E$ is:

\begin{equation}\label{EquationVISFMState}
x=[r,~\phi,~v,~\alpha,~\theta,~A_y]^T
\end{equation}
The unknown inputs of $\E$ are $A_x$ and $\dt{A}_y$.
By using (\ref{EquationVISFMDynamics}), we easily obtain the dynamics of the new state and, by comparing with (\ref{EquationSystemDefinitionUIO}), we obtain:

\begin{equation}\label{EquationVSFMg}
g^0=\left[\begin{array}{c}
v \cos(\alpha-\phi)\\
\frac{v}{r} \sin(\alpha-\phi)\\
 A_y\sin(\alpha-\theta)   \\
\frac{A_y}{v} \cos(\alpha-\theta)\\
 0  \\
 0  \\
\end{array}
\right],~
g^1=\left[\begin{array}{c}
 0 \\
 0  \\
 \cos(\alpha-\theta)   \\
-\frac{1}{v} \sin(\alpha-\theta)\\
 0  \\
 0  \\
\end{array}
\right]
\end{equation}
\[
g^2=[0, 0, 0, 0, 0, 1]^T, f^1=[0, 0, 0, 0, 1, 0]^T.
\]

Algorithm 9 in \cite{IF22} also provides
the unknown input degree of reconstructability of $\E$, which is $m=2$,
\begin{equation}\label{EquationVSFMht}
\widetilde{h}_1=\frac{v}{r} \sin(\alpha-\phi),~~
\widetilde{h}_2=A_y\frac{\cos(\phi-\theta)}{ r},
\end{equation}
and the observability codistribution of $\E$, which is $
\OBS=\textnormal{span}\left\{
d\phi - d\theta,~~
d\widetilde{h}_1,~~
d\widetilde{h}_2,~~
dh_3
\right\}$, with $h_3=A_y\frac{\sin(\phi-\theta)}{ r}$.
We have $m=2=m_w$, meaning that $\E$ is canonic with respect to its unknown inputs.

\subsection{Identifiability analysis}\label{SubSectionVISFMIdentifiability}

Our system is characterized by two time varying parameters ($A_x$ and $\dt{A}_y$). In addition it is canonic with respect to them. As a result, all the identifiability properties are obtained by using the results of Theorem \ref{TheoremIdentifiabilityGeneral}. 
First of all, we need to compute the orthogonal distribution $(\OBS)^{\bot}$. We obtain:

\begin{equation}\label{EquationVISFMOBSOrtho}
(\OBS)^{\bot}=\textnormal{span}\{
[0~1~0~1~1 ~0]^T,~~
[r~0~v~0~0~A_y]^T
\}.
\end{equation}
Note that, from the expression of its second generator, we immediately obtain that $dA_y\notin\OBS$, meaning that the original second time varying parameter $A_y(t)$ is not locally identifiable.
To check the local identifiability  of $A_x$ and $\dt{A}_y$ we use the condition in (\ref{EquationConditionGeneral}).
We compute the $(1,1)$-tensor $\mu$ in (\ref{EquationConditionMu}) by using (\ref{EquationVSFMg}) and (\ref{EquationVSFMht}). We obtain for its inverse:

\[
\nu^1_1= -\frac{ r}{\sin(\phi - \theta)},~~
\nu^1_2=~\nu^2_1=0, ~~
\nu^2_2=\frac{ r}{\cos(\phi - \theta)}.
\]

Then, for each generator of $(\OBS)^{\bot}$ in (\ref{EquationVISFMOBSOrtho}), we compute the coefficients $\xi^\alpha_i$, for $\alpha=0,1,2$, and $i=1,2$, by using (\ref{EquationConditionCoeffXi}), (\ref{EquationVSFMg}), (\ref{EquationVSFMht}), and (\ref{EquationVISFMOBSOrtho} ). We obtain that, for the first generator, all of them vanish. As a result, by using the first generator, the condition in (\ref{EquationConditionGeneral}) is honoured.
Regarding the second generator, we obtain:

\[
\xi^0_1=\xi^0_2=\xi^1_2=\xi^2_1=0,~
\xi^1_1=\frac{\sin(\phi-\theta)}{r},~
\xi^2_2=-\frac{\cos(\phi-\theta)}{r}
\]
By substituting the above $\nu^i_j$ and $\xi^\alpha_i$ in (\ref{EquationConditionGeneral}), we have the following result. For the first parameter $A_x(t)$ (i.e., for $j=1$ in (\ref{EquationConditionGeneral})), we obtain a non vanishing result only for $\alpha=1$. Specifically, we obtain $-1$.
~For the second parameter $\dt{A}_y(t)$ (i.e., for $j=2$ in (\ref{EquationConditionGeneral})), we obtain a non vanishing result only for $\alpha=2$. Specifically, we obtain again $-1$.

We conclude that both the time varying parameters are not locally identifiable.

Finally, we use the second part of Theorem \ref{TheoremIdentifiabilityGeneral} to determine infinitely many values of these parameters indistinguishable from their true values. From the above, we know that we must use only the second generator to build the system of differential equations in (\ref{EquationConditionDiffEqSystem}). The components of $\chi$ in (\ref{EquationConditionChij}) are:
\[
\chi_1=A_x,~~\chi_2=\dot{A}_y
\]
and the system of differential equations in (\ref{EquationConditionDiffEqSystem}) is:
\begin{equation}\label{EquationVISFMDiffEqSystem}
\left\{\begin{array}{ll}
  \frac{dr'}{d\tau} &=  r'\\
  \frac{d\phi'}{d\tau} &=  0\\
  \frac{dv'}{d\tau} &=  v'\\
  \frac{d\alpha'}{d\tau} &=  0\\
  \frac{d\theta'}{d\tau} &=  0\\
  \frac{dA'_y}{d\tau} &=  A'_y\\
  \frac{dA'_x}{d\tau} &=  A'_x\\
  \frac{d\dt{A}'_y}{d\tau} &=  \dt{A}'_y\\
r'(t,~0)&=r(t), ~\phi'(t,~0)= \phi(t)\\
v'(t,~0)&=v(t), ~\alpha'(t,~0)= \alpha(t)\\
\theta'(t,~0)&=\theta(t), ~A'_y(t,~0)= A_y(t)\\
A'_x(t,~0)&=A_x(t), ~\dt{A}'_y(t,~0)= \dt{A}_y(t)\\
\end{array}\right.
\end{equation}
Its solution is trivially obtained and it is:
\begin{equation}\label{EquationVISFMIndistinguishableStatesAndUI}
r'(t,~\tau)=e^\tau r(t), ~\phi'(t,~\tau)= \phi(t),~v'(t,~\tau)=e^\tau v(t)
\end{equation}
\[
\alpha'(t,~\tau)= \alpha(t), ~\theta'(t,~\tau)=\theta(t),~A'_y(t,~\tau)= e^\tau A_y(t)
\]
\[
A'_x(t,~\tau)= e^\tau A_x(t), ~\dt{A}'_y(t,~\tau)= e^\tau\dt{A}_y(t)
\]

In other words, starting from the true values of the state and the time varying parameters, we can determine infinitely many values indistinguishable from them by performing a scale transform (the scale is precisely $e^\tau$, for any $\tau\in\R$).
This result is not surprising. All the information is provided by the measurements that only consist of angular measurements (the visual sensor only provides the angle $\beta$ in Fig. \ref{FigVISFM} and the gyroscope only provides the angular speed). As a result, the system has no source of metric information and the determination of any metric quantity is only possible up to a scale.

\section{Conclusion}\label{SectionConclusion}
In this paper we provided the following three main contributions:

\begin{itemize}

\item {\bf First contribution:} Introduction of the general analytical condition that fully characterizes the local identifiability of all the unknown parameters of any ODE model.

\item {\bf Second contribution:} Introduction of the continuous transformations that allow us to determine infinitely many values, for any unidentifiable parameter (constant or time varying), that are indistinguishable from its true value.

\item {\bf Third contribution:}
Introduction of new important results about the identifiability 
of a very popular HIV model and a Covid-19 model. The results on the HIV model highlight
a serious error in the state of the art.


\end{itemize}

\subsubsection{First contribution}

The analytical condition was provided in Section \ref{SectionCondition}. 
It can be applied to any system, regardless of its complexity and type of nonlinearity and in the presence of time varying parameters. In particular, it can be applied to the very general model in  (\ref{EquationSystemDefinitionUIOGeneral}), under Assumptions \ref{Assumptionfh}-\ref{AssumptionU}.
Its application requires no inventiveness from the user. It is sufficient to carry out very basic and automatic computation (derivatives and matrix ranks).

\subsubsection{Second contribution}

When one or more parameters (constant and/or time varying) are unidentifiable, there are more values of them that reproduce exactly the same outputs (and also agree with the same known inputs, when present).
Section \ref{SectionCondition} introduced the mathematical tool that quantitatively allows us to determine these indistinguishable values. This tool consists of a system of first order ordinary differential equations, which is given in (\ref{EquationConditionDiffEqSystem}).

\subsubsection{Third contribution}

Sections \ref{SectionHIV} and \ref{SectionCovid} provided a detailed study of the identifiability properties of two viral ODE models. In particular, the ODE model investigated in Section \ref{SectionHIV} is very popular and regards the HIV dynamics. 
%
%
%
%
%
~For this model, the system of differential equations in (\ref{EquationConditionDiffEqSystem}) becomes the system in (\ref{EquationHIVDiffEqSystem}). This system is complex and it is not surprising that, in the state of the art, no value of the  time varying parameter, different from the true one but indistinguishable from it, was detected, and the model was classified {\it identifiable}. In other words, the solution of (\ref{EquationHIVDiffEqSystem}), although obtained analytically, is complex and not possible to be determined by following an intuitive reasoning, as for the visual inertial sensor fusion problem, studied in Section \ref{SectionVISFM}. In this last case, the system of differential equations in (\ref{EquationConditionDiffEqSystem}) becomes the system in (\ref{EquationVISFMDiffEqSystem}), which is trivial. Its solution is (\ref{EquationVISFMIndistinguishableStatesAndUI}) and could have been obtained with intuitive reasoning (basically, the physical meaning of this solution is that we cannot reconstruct the absolute scale, which was intuitive).

\vskip.2cm

We conclude by emphasizing the generality of the results introduced by this paper. They provide, automatically, the local identifiability of all the parameters (constant and time varying) of any nonlinear system, independently of its complexity and type of nonlinearity. 
This paper provided a first application of all these theoretical results on two viral models and several important new properties (and also an error in the state of the art) were automatically discovered.
Due to the aforementioned generality, we strongly believe that our theoretical results can be very useful in many other scientific domains, basically to investigate any system that can be described by a an ODE model with time varying parameters. This occurs in many scientific domains, ranging from the natural and applied sciences up to the social sciences.

\appendix

This appendix provides the proof of Theorem \ref{TheoremIdentifiabilityCan} (Section \ref{SubSectionAppendixProofTheoremIdentifiabilityCan}), of Theorem \ref{TheoremIdentifiabilityGeneral} (Section \ref{SubSectionAppendixProofTheoremIdentifiabilityGeneral}) and of Theorem \ref{TheoremIdentifiabilityObservable} (Section \ref{SubSectionAppendixProofTheoremIdentifiabilityObservable}). We start by characterizing the concept of identifiability and local identifiability (Section \ref{SubSectionAppendixDefinitionIdentifiability}) in the presence of time varying parameters and when the model has an explicit time dependence.


\subsection{Definition of Identifiability}\label{SubSectionAppendixDefinitionIdentifiability}

We introduce the definition of identifiability for the model in (\ref{EquationSystemDefinitionUIOGeneral}). 
Clearly, this definition also holds for the model in (\ref{EquationSystemDefinitionUIO}), which is equivalent under Assumption \ref{AssumptionU}. We only emphasize that, in this second case, the definition of the $\theta$-state given below (Equation (\ref{EquationAppendixThetaState})) does not need to include the constant parameters because they are already in $x$.

The concept of identifiability will be based on the concept of indistinguishability. 
We must account for (i) the explicit time dependence of the model, and (ii) the presence of time varying parameters. As the time varying parameters  act on the system dynamics as unknown inputs, we use the same concepts which were adopted to introduce the concept of indistinguishability in the presence of unknown inputs (e.g., see Section II in \cite{SARAFRAZI} or Definition 6.15 in \cite{SIAMbook}). In addition, 
we account for the explicit time dependence of the model in (\ref{EquationSystemDefinitionUIOGeneral}) by proceeding as in \cite{TAC22}. On the other hand, in contrast with these previous works, we are here interested in the identifiability of the parameters instead of the state observability. However, the observability of the state and the identifiability of the parameters are interconnected. To account for this, we start by introducing an extended state that includes the state and all the parameters. We use the symbol $\theta$ to denote all the parameters (i.e., $\theta=\left[Q^T,~W^T(t)\right]^T$).
We introduce the $\theta$-state, defined as follows:

\begin{equation}\label{EquationAppendixThetaState}
S=\left[X^T,~\theta^T\right]^T=\left[X^T,~Q^T,~W^T\right]^T.
\end{equation}
Note that, for the model in (\ref{EquationSystemDefinitionUIO}) we have $S=\left[x^T,~w_1,\ldots,w_{m_w}\right]^T$.

We denote by $X(t; ~t_0, \overline{X}, \overline{Q}, \overline{W})$ the state at the time $t\ge t_0$ that satisfies the dynamics in (\ref{EquationSystemDefinitionUIOGeneral}) when the model parameters are 
$Q=\overline{Q}$, $W(t)=\overline{W}(t)$, and the state at tine $t_0$ is $\overline{X}$.

Let $\mathcal{W}$
be the Banach function space of all possible first order time derivative of the unknown input functions ($\dt{W}$). 

We define indistinguishability as follows:

\begin{df}[Indistinguishable $\boldsymbol{\theta}$-states]\label{DefinitionIndistinguishableThetaStates}
Given the system in (\ref{EquationSystemDefinitionUIOGeneral}), two $\theta$-states, $S_a=\left[X_a^T,~Q_a^T,W_a^T\right]^T$, $S_b=\left[X_b^T,~Q_b^T,W_b^T\right]^T$, are indistinguishable at $t_0\in\I$
if, for any time assignment $U(t)$, $t\in\I\bigcap[t_0,\infty)$, there exists a nonshy\footnote{See \cite{Hunt92} for a definition of shy set.} subset $\mathcal{W}_a\subseteq\mathcal{W}$ such that, for any 
unknown input $W_A(t)$ with $\dt{W}_A\in\mathcal{W}_a$ and $W_A(t_0)=W_a$, there exists an unknown input $W_B(t)$ with $W_B(t_0)=W_b$ such that, $\forall t\in\I\bigcap[t_0,\infty)$:
\[
h(X(t; ~t_0,X_a, Q_a, W_A), ~t, ~U(t), ~Q_a, ~W_A(t))=
\]
\[h(X(t; ~t_0,X_b, Q_b, W_B), ~t, ~U(t), ~Q_b, ~W_B(t)).
\]
\end{df}
%
The above definition takes into account that we do not know the time behaviour of the unknown input $W(t)$. In particular, Definition \ref{DefinitionIndistinguishableThetaStates} does not exclude the possibility that the outputs produced by $S_a$ and $S_b$ can differ even for many time profiles $W_A(t)$ and $W_B(t)$. However, Definition \ref{DefinitionIndistinguishableThetaStates} establishes that the outputs coincide for all the profiles that belong to a nonshy set. In other words, there exists a non-vanishing probability that the outputs coincide\footnote{For finite-dimensional spaces, the probability is defined by using the concept of measure.
Unfortunately, the concept of measure in spaces with infinite dimensions is not trivial. 
In particular, there is no analogue of Lebesgue measure on an infinite-dimensional Banach space, such as our function space $\mathcal{W}$. One possibility, which is frequently adopted, is to use the concept of prevalent and shy sets \cite{Hunt92}. The probability that a given $\dt{W}$ belongs to a shy subset of $\mathcal{W}$, is $0$. The probability that a given $\dt{W}$ belongs to a prevalent subset of $\mathcal{W}$, is $1$. Finally, the probability that a given $\dt{W}$ belongs to a nonshy subset of $\mathcal{W}$, is strictly larger than $0$.}. This will make conservative the following definitions of identifiability and local identifiability (Definitions \ref{DefinitionIdentifiability} and \ref{DefinitionLocalIdentifiability}), because they are based on Definition \ref{DefinitionIndistinguishableThetaStates}.

\begin{df}[Identifiability]\label{DefinitionIdentifiability}
The system in (\ref{EquationSystemDefinitionUIOGeneral}) is identifiable at a given $X\in\M_m$, $Q\in\M_q$, $W\in\M_{m_w}$, and $t\in\I$ if all the $\theta$-states $\left[\overline{X}^T,~\overline{Q}^T,\overline{W}^T\right]^T$ indistinguishable from $\left[X^T,~Q^T,~W^T\right]^T$ at $t$, are such that:
\begin{equation}\label{EquationDefinitionIdentifiability}
\overline{Q}=Q,~~\textnormal{and}~~\overline{W}=W.
\end{equation} 
It is identifiable on a subset $\mathcal{G}\subseteq\M_m\times\M_q\times\M_{m_w}\times\I$, if it is identifiable at any $(X,~Q,~W~,t)\in\mathcal{G}$.
\end{df}

\begin{df}[Local identifiability]\label{DefinitionLocalIdentifiability}
The system in (\ref{EquationSystemDefinitionUIOGeneral}) is locally identifiable at a given $X\in\M_m$, $Q\in\M_q$, $W\in\M_{m_w}$, and $t\in\I$ if 
there exists an open neighbour $\mathcal{B}$ of $\left[X^T,~Q^T,~W^T\right]^T$ such that
all the $\theta$-states $\left[\overline{X}^T,~\overline{Q}^T,\overline{W}^T\right]^T\in\mathcal{B}$ indistinguishable from $\left[X^T,~Q^T,~W^T\right]^T$ at $t$, are such that the same equalities in (\ref{EquationDefinitionIdentifiability}) hold.
It is locally identifiable on a subset $\mathcal{G}\subseteq\M_m\times\M_q\times\M_{m_w}\times\I$, if it is locally identifiable at any $(X,~Q,~W~,t)\in\mathcal{G}$.
\end{df}

The two above definitions can be easily extended in order to characterize the identifiability (or the local identifiability) of a single parameter. Let us denote by $\theta_l$, $l=1,\ldots,q+m_w$, one of the parameters of our model in (\ref{EquationSystemDefinitionUIOGeneral}). In particular, $\theta_1,\ldots,\theta_q$ are the constant parameters $Q_1,\ldots,Q_q$, and $\theta_{q+1},\ldots,\theta_{q+m_w}$ are the time varying parameters $W_1,\ldots,W_{m_w}$. The definition of the identifiability (or the local identifiability) of the parameter $\theta_l$ is obtained by simply replacing the equalities in (\ref{EquationDefinitionIdentifiability}) with the single equality $\overline{\theta}_l=\theta_l$.


\subsection{Proof of Theorem \ref{TheoremIdentifiabilityCan}}\label{SubSectionAppendixProofTheoremIdentifiabilityCan}


We consider the system $\E$ provided by Algorithm 9 in \cite{IF22}, with $m$ its unknown input degree of reconstructability, and $\OBS$ its observability codistribution. $\E$ is characterized by (\ref{EquationSystemDefinitionUIO}) and may differ from the original system, as explained in Section \ref{SubSectionConditionReminders}.

We start by proving the following property:

\begin{lm}\label{LemmaUIRecLieZeroLast}
If $m<m_w$ we have:

\[
\Li_{g^{m+j}}\gamma=0,~~j=1,\ldots,m_w-m,
\]
where $g^1,\ldots,g^{m_w}$ are the vector fields that characterize the dynamics of $\E$ and $\gamma$ is any scalar function of the state (of $\E$) such that $d\gamma\in\OBS$.
\end{lm}

\proof{We remind the reader that, when the system is not canonic with respect to its unknown inputs, the observability codistribution is the limit codistribution of the series $\OBS^\infty_k$ defined in Section 6.2 of \cite{IF22}, which are computed by Algorithm 8 in \cite{IF22}.
From Proposition 2 in \cite{IF22} we know that the series of codistributions $\OBS_k^\infty$ converges and the limit codistribution is precisely $\OBS$.
Hence, the codistribution $\OBS$ is invariant under the Lie derivatives along the vector fields $^m\widehat{g}^\alpha_{d\infty}$, which are defined by Equation (31) in \cite{IF22}, for any $\alpha=0,1,\ldots,m$.
As a result, for any function $\gamma$, such that $d\gamma\in\OBS$, we obtain that $\Li_{^m\widehat{g}^\alpha_{d\infty}}\gamma$
is independent of $w_{m+j}$, for any $j=1,\ldots,m-m_w$ and for any $\alpha=0,1,\ldots,m$ (note that, for a time varying  system, and for $\alpha=0$, this result actually regards $\dt{\Li}_{^m\widehat{g}^0_{d\infty}}\gamma$, with $\dt{\Li}_{^m\widehat{g}^0_{d\infty}}:=\Li_{^m\widehat{g}^0_{d\infty}}+\frac{\partial}{\partial t}$).
Let us consider the above function when $\alpha=0$. 
From Equation (31) in \cite{IF22}, and by knowing that $^m\nu^0_0=1$, we have:
 
 \[
 ^m\widehat{g}^0_{d\infty}=g^0_{d\infty} + \sum_{j=1}^d~^m\nu^0_jg^j_{d\infty},
 \]

with $d=m_w-m$.
From Equations (27) and (28) in \cite{IF22}, we remark that, in the above expression, only
$g^0_{d\infty}$ depends on $w_{m+j}$ ($j=1,\ldots,m_w-m$). 
On the other hand:

\[
\Li_{^m\widehat{g}^0_{d\infty}}\gamma=\Li_{g^0_{d\infty}}\gamma +  \sum_{j=1}^d~^m\nu^0_j \Li_{g^j_{d\infty}}\gamma
\]

and, in this function, only the first term,
$\Li_{g^0_{d\infty}}\gamma$, can depend on $w_{m+j}$.

From Equation (27) in \cite{IF22}, we have:
\[
\Li_{g^0_{d\infty}}\gamma=\Li_{g^0}\gamma+\sum_{j=1}^d (\Li_{g^{m+j}}\gamma)~ w_{m+j}
\]

Since this function must be independent of  $w_{m+j}$, for any $j=1,\ldots,m-m_w$, we obtain:
\[
\Li_{g^{m+j}}\gamma=0
\]
for any $j=1,\ldots,m-m_w$ and for any function $\gamma$ such that $d\gamma\in\OBS$.
$\blacktriangleleft$}

\vskip.1cm

Now, we are ready to prove Theorem \ref{TheoremIdentifiabilityCan}.

\vskip.1cm

\proof{
The definition of local identifiability of a given parameter provided in Appendix \ref{SubSectionAppendixDefinitionIdentifiability}, coincides with the definition given for the observability in the presence of unknown inputs and for time varying systems (e.g., see Definition 6.15 in \cite{SIAMbook} or Section II in \cite{SARAFRAZI}), when applied to the corresponding component of the $\theta$-state.
This allows us to study the identifiability of a given parameter by studying the observability codistribution of the $\theta$-state. 
Let us denote this codistribution by $\OBS^S$. The state is (we are referring to the system $\E$ which is characterized by (\ref{EquationSystemDefinitionUIO})):

\begin{equation}\label{EquationAppendixProofTheoremIdCanState}
S = [x^T,~ w_1,\ldots,w_{m_w}]^T
\end{equation}
where $x$ is the state of $\E$ and it already includes all the constant parameters.
The local identifiability of a given parameter will be obtained by simply verifying whether or not its differential belongs to $\OBS^S$. 

Let us compute $\OBS^S$. We know $\OBS$, which is the observability codistribution of $\E$.

Starting from the expression of the vector fields that characterize the dynamics of $x$, it is immediate to obtain the expression of the vector fields that characterize the dynamics of $S$. In particular, the drift, which will be denoted by $G^0$, has the following expression:

\begin{equation}\label{EquationAppendixProofTheoremIdCanG0}
G^0=
\left[
\begin{array}{c}
g^0+\sum_{j=1}^{m_w}g^jw_j\\
0_{m_w}\\
\end{array}
\right]
\end{equation}

where $0_{m_w}$ is the zero $m_w$-column vector.

We know that, for any scalar function $\gamma(x)$ such that $d\gamma\in\OBS$, we have (only for the first order Lie derivative):

\[
d\dt{\Li}_{G^0}\gamma\in\OBS^S.
\]
where $\dt{\Li}_{G^0}=\Li_{G^0}+\frac{\partial}{\partial t}$\footnote{In the time invariant case, we have $d\Li_{G^0}\gamma\in\OBS^S$. This property is well known and widely used in the past (e.g., in all the procedures which are based on a state extension introduced in \cite{Belo10,MED15,Villa19b,Maes19,Villa16b}). In the time varying case, we need to add the time derivative operator (see Algorithm 2, Theorem 2, and Equation (10) in \cite{TAC22}).}.


We denote the rank of $\OBS$ by $n_o$ ($\le n$, where $n$ is the dimension of $x$). Algorithm 9 in \cite{IF22} provides the $n_o$ generators of $\OBS$. Let us denote them by 
$\gamma_1(x),\ldots,\gamma_{n_o}(x)$.
We have:

\[
\OBS=\textnormal{span}\left\{
d\gamma_1,\ldots,d\gamma_{n_o}
\right\}.
\]

Algorithm 9 in \cite{IF22} also provides the functions $\widetilde{h}_1,\ldots,\widetilde{h}_m$, which are $m$ scalar fields such that the reconstructability matrix has full rank and whose differentials belong to $\OBS$, i.e., 
$d\widetilde{h}_i\in\OBS$, $i=1,\ldots,m$. As a result:

\[
d\dt{\Li}_{G^0}\widetilde{h}_i\in\OBS^S, ~~~i=1,\ldots,m.
\]

In addition, as the the reconstructability matrix from them has full rank, the covectors $d\dt{\Li}_{G^0}\widetilde{h}_1,\ldots,d\dt{\Li}_{G^0}\widetilde{h}_m$ are independent and they are also independent from the generators of $\OBS$ (because the latter are independent of the unknown inputs).
As a result, the $n_o+m$ covectors $d\gamma_1,\ldots,d\gamma_{n_o}, d\dt{\Li}_{G^0}\widetilde{h}_1,\ldots,d\dt{\Li}_{G^0}\widetilde{h}_m$ are independent.


From Lemma \ref{LemmaUIRecLieZeroLast}, and from the structure of $G^0$ in (\ref{EquationAppendixProofTheoremIdCanG0}), we obtain that, for any scalar function $\gamma(x)$ such that $d\gamma\in\OBS$, $\dt{\Li}_{G^0}\gamma$  is independent of the last $m_w-m$ unknown inputs (i.e., the last $m_w-m$ entries of $S$). As a result, the rank of $\OBS^S$ cannot exceed $n_o+m$.
Therefore:

\begin{equation}\label{EquationAppendixProofTheoremIdCanOBS}
\OBS^S=\textnormal{span}\left\{
d\gamma_1,\ldots,d\gamma_{n_o}, ~d\dt{\Li}_{G^0}\widetilde{h}_1,\ldots,d\dt{\Li}_{G^0}\widetilde{h}_m
\right\}.
\end{equation}


Let us study the structure of the null space of $\OBS^S$.
As all the $n_o+m$ functions $\gamma_1,\ldots,\gamma_{n_o},\dt{\Li}_{G^0}\widetilde{h}_1,\ldots,\dt{\Li}_{G^0}\widetilde{h}_m$
 are independent of the last $m_w-m$ unknown inputs, any vector with dimension $n+m_w$, with the first $n+m_w-m$ entries equal to zero, certainly belongs to the null space of $\OBS^S$. Therefore, the differential $dw_{m+j}\notin\OBS^S$.
This proves the first part of the Theorem.

To prove the second part, we exploit the following property. 
Let us denote by $\Delta^S$ the distribution orthogonal to the observability codistribution (i.e., $\Delta^S:=\left(\OBS^S\right)^\bot$). $\Delta^S$ is involutive ($\OBS^S$ is an integrable codistribution by construction).
Let us denote by $S(t)$ the true $\theta$-state at time $t$. For any vector field $\xi^S\in\Delta^S$, let us consider the following system of differential equations in $\tau$ (where $t$ is a fixed parameter):

\begin{equation}\label{EquationAppendixDiffEqSystemPROPRIETA}
\left\{\begin{array}{ll}
  \frac{dS'}{d\tau} &=   \xi^S(S'(t,~\tau))\\
  S'(t,~0) &=S(t),\\
\end{array}\right.
\end{equation}
Let us assume that the solution of (\ref{EquationAppendixDiffEqSystemPROPRIETA}) exists for any $t\in\I$ and for any $\tau\in\mathcal{U}$, with $\mathcal{U}$ an open interval of $\R$ that includes $\tau=0$. Then, the solution of (\ref{EquationAppendixDiffEqSystemPROPRIETA}) at any $\tau\in\mathcal{U}$, i.e., $S'(t,~\tau)$, is indistinguishable from $S(t)$\footnote{This can be easily proved by showing that any scalar function $\Gamma(S)$ of the $\theta$-state such that its differential belongs to the observability codistribution (i.e., $d\Gamma\in\OBS^S$) satisfies $\frac{\partial\Gamma(S'(t,~\tau))}{\partial\tau}=0$. The validity of this equality is an immediate consequence of the orthogonality between $d\Gamma$ and $\xi^S$.}.

Let us consider one of the last $m_w-m$ unknown inputs, i.e., the unknown input $w_{m+k}$, for a given $k=1,\ldots,m_w-m$.

In the first part of this proof, we proved that any vector with dimension $n+m_w$, with the first $n+m_w-m$ entries equal to zero, certainly belongs to the null space of $\OBS^S$. As a result, the following vector field certainly belongs to $\Delta^S$:

\[
\xi^S_k=\left[
\begin{array}{l}
0_n\\
0_m\\
e_k\\
\end{array}
\right]
\]
with $0_n$ the zero $n$-column vector, $0_m$ the zero $m$-column vector, and $e_k$ the $(m_w-m)$-column vector, with component $k^{th}$ 1 and zero elsewhere.
By replacing $\xi^S$ in (\ref{EquationAppendixDiffEqSystemPROPRIETA})
with the above $\xi^S_k$, by a direct integration of (\ref{EquationAppendixDiffEqSystemPROPRIETA}), we immediately obtain that 
$w_{m+k}'(t,~\tau)=w_{m+k}(t)+\tau$ is indistinguishable from $w_{m+k}(t)$, $\forall\tau\in\R$.
$\blacktriangleleft$}

\subsection{Proof of Theorem \ref{TheoremIdentifiabilityGeneral}}\label{SubSectionAppendixProofTheoremIdentifiabilityGeneral}

\proof{
We proceed as in the proof of Theorem \ref{TheoremIdentifiabilityCan}.
In this case, we first prove the second part of the theorem.
It is sufficient to prove that, $\forall\xi\in\OBS^\bot$, the vector field:

\begin{equation}\label{EquationAppendixSymmetryS}
\xi^S:=
\left[
\begin{array}{c}
\xi\\
\chi\\
\end{array}
\right],
\end{equation}
with $\chi$ defined in (\ref{EquationConditionChij}),
belongs to the null space of $\OBS^S$ in (\ref{EquationAppendixProofTheoremIdCanOBS}).

Let us consider the first $n_o$ generators of $\OBS^S$. As $\xi\in\OBS^T$, we have:

\[
\frac{\partial\gamma_i}{\partial x}\cdot\xi=0, ~~~i=1,\ldots,n_o.
\]

Hence,

\[
\frac{\partial\gamma_i}{\partial S}\cdot\xi^S=\frac{\partial\gamma_i}{\partial x}\cdot\xi+\frac{\partial\gamma_i}{\partial w}\cdot~\chi=
\]
\[
0+
\underbrace{\left[0,\ldots,0\right]}_{m_w}
\cdot~\chi=0, ~~~i=1,\ldots,n_o.
\]
as $\gamma_i$ is independent of $w$ (it only depends on $x$).

Let us consider the remaining $m$ generators. We have:

\begin{equation}\label{EquationAppendixProofTheoremGeneralLiG0ht}
\frac{\partial\dt{\Li}_{G^0}\widetilde{h}_l}{\partial S}\cdot\xi^S=
\frac{\partial \dt{\Li}_{G^0}\widetilde{h}_l}{\partial x}\cdot\xi+
 \sum_{k=1}^m\frac{\partial \dt{\Li}_{G^0}\widetilde{h}_l}{\partial w_k}\chi_k
\end{equation}

On the other hand,

\begin{equation}\label{EquationAppendixProofTheoremGeneralLiG0ht1}
\dt{\Li}_{G^0}\widetilde{h}_l=\dt{\Li}_{g^0}\widetilde{h}_l + \sum_{j=1}^m (\Li_{g^j}\widetilde{h}_l) w_j,~~~l=1,\ldots,m
\end{equation}

By using (\ref{EquationAppendixProofTheoremGeneralLiG0ht1}) in (\ref{EquationAppendixProofTheoremGeneralLiG0ht}) and by using (\ref{EquationConditionMu}) and (\ref{EquationConditionCoeffXi}) we obtain:

\[
\frac{\partial\dt{\Li}_{G^0}\widetilde{h}_l}{\partial S}\cdot\xi^S=
\xi^0_l+ \sum_{j=1}^m \xi^j_l w_j +
\sum_{k=1}^m \Li_{g^k}\widetilde{h}_l \chi_k=
\]
\[
\xi^0_l+ \sum_{j=1}^m \xi^j_l w_j +
\sum_{k=1}^m ~\mu^k_l ~\chi_k.
\]

Now, by using the expression of $\chi_k$ in (\ref{EquationConditionChij}) and that $\mu$ is the inverse of $\nu$ (i.e., $\sum_{k=1}^m~\mu^k_l ~\nu^i_k=\delta^i_l$), we finally obtain:

\[
\frac{\partial\dt{\Li}_{G^0}\widetilde{h}_l}{\partial S}\cdot\xi^S=0.
\]
Hence, we proved that $\xi^S\in\Delta^S$. By using this in (\ref{EquationAppendixDiffEqSystemPROPRIETA}), we obtain the proof of the second part of the theorem.

Let us prove the first part. As we aforementioned, checking the local identifiability of a given parameter is equivalent to check whether or not its differential belongs to $\OBS^S$. But this is equivalent to check if the component of $\chi$ that corresponds to this parameter vanishes for all the vector fields that belong to the null space of $\OBS^S$.
In other words, regarding the $j^{th}$ time varying parameter, its local identifiability is equivalent to the condition $\chi_j=0$,
for any $\xi\in\OBS^\bot$, with $\chi_j$ given in (\ref{EquationConditionChij}). As all the coefficients $\xi^\alpha_i$ that appear in (\ref{EquationConditionChij}) are independent of $w$ (they only depend on $x$), $\chi_j=0$ if and only if we have the condition in (\ref{EquationConditionGeneral}), for any $\alpha=0,1,\ldots,m$.
$\blacktriangleleft$}

\subsection{Proof of Theorem \ref{TheoremIdentifiabilityObservable}}\label{SubSectionAppendixProofTheoremIdentifiabilityObservable}

\proof{
Based on the result stated by Theorem 4 in \cite{IF22}, it remains to prove the viceversa, i.e., if the unknown inputs are locally identifiable then the system is canonic with respect to its unknown inputs. 

We proceed by contradiction. We assume that the system is not canonic with respect to its unknown inputs.

As the state is observable, all its components belong to the observation space and we have $\OBS=\textnormal{span}\left\{
d x_1, \ldots, d x_n\right\}$, at least locally.
As we assumed that the system is not canonic, $\uideg\left(\OBS\right)<m_w$. Therefore, the matrix $\RM\left( x_1,\ldots, x_n\right)$ has rank smaller than $m_w$. On the other hand, by an explicit computation we obtain:

\[
\RM\left( x_1,\ldots, x_n\right)
=
\left[\begin{array}{cccc}
~[g^1]_1 & [g^2]_1 & \ldots & [g^{m_w}]_1 \\
~[g^1]_2 & [g^2]_2 & \ldots & [g^{m_w}]_2 \\
\ldots &\ldots &\ldots &\ldots \\
~[g^1]_n & [g^2]_n & \ldots & [g^{m_w}]_n \\
\end{array}
\right]=
\]
\[
\left[\begin{array}{cccc}
g^1 & g^2 & \ldots & g^{m_w} \\
\end{array}
\right],
\]
where $[g^j]_i$ is the $i^{th}$ component of the vector field $g^j$ ($j=1,\ldots,m_w$ and $i=1,\ldots,n$).
Therefore, as rank$\left(\RM\left( x_1,\ldots, x_n\right)\right)<m_w$, the vectors $g^1 ~g^2 ~\ldots ~ g^{m_w}$ are linearly dependent. 
We denote by $\widehat{n}$ a non trivial vector that belongs to the null space of the above matrix. 
Let us consider the dynamics:

\[
\dot{x} =   g^0+\sum_{k=1}^{m_u}f^k  u_k +  \sum_{j=1}^{m_w}g^j  w_j.
\]

For any $w=[w_1,~w_2,\ldots,w_{m_w}]$, we obtain the same dynamics by transforming the UI as follows:

\[
w\rightarrow w'=w+\widehat{n}
\]

As result, $w$ cannot be distinguished from $w'$.
$\blacktriangleleft$}

%
%
%
%
%

\end{document}